\documentclass[authoryear,11pt,preprint]{elsarticle}

\usepackage{color}
\usepackage{lineno}
\usepackage{pdflscape}
\usepackage{amsmath}
\usepackage{booktabs}
\usepackage[hyphens]{url}
\usepackage{hyperref}
\usepackage{multirow}
\usepackage{addlines}
\usepackage{adjustbox}
\usepackage{natbib} 
\usepackage[english]{babel}
\usepackage[dvipsnames]{xcolor}
\usepackage{soul}
\usepackage{adjustbox}
\usepackage{verbatim}
\usepackage[a4paper, top=2.5cm, bottom=2.5cm, left=2.5cm, right=2.5cm]{geometry}

\usepackage{mathtools}
\usepackage{tikz}
\usepackage{pgfplots}
\usepackage{enumitem}
\usepackage{calc}

\definecolor{darkblue}{rgb}{0,0.2,0.4}
\usepackage{graphicx}
\usepackage{subcaption}
\usepackage{amssymb}
\usepackage{amsthm}
\usepackage{etoolbox}
\usepackage{algorithm}
\usepackage{algpseudocode}

\usepackage{float} 
\usepackage{algpseudocode} 
\usepackage{algorithm} 

\newfloat{subroutine}{htbp}{loa}
\floatname{subroutine}{Subroutine}

\makeatletter
\newcommand*{\algrule}[1][\algorithmicindent]{%
  \hspace*{.2em}
  \vrule 
  \hspace*{\dimexpr#1-.2em-.4pt}%
}

\newcommand{\StatePar}[1]{%
  \State\parbox[t]{\dimexpr\linewidth-\ALG@thistlm}{\strut #1\strut}%
}
\renewcommand{\ALG@beginalgorithmic}{\offinterlineskip}

\newcount\ALG@printindent@tempcnta
\def\ALG@printindent{%
  \ifnum \theALG@nested > 0
    \ifx\ALG@text\ALG@x@notext
    \else
      \unskip
      \ALG@printindent@tempcnta=1
      \loop
        \algrule[\csname ALG@ind@\the\ALG@printindent@tempcnta\endcsname]%
        \advance \ALG@printindent@tempcnta 1
        \ifnum \ALG@printindent@tempcnta<\numexpr\theALG@nested+1\relax
      \repeat
        \fi
    \fi
}
\patchcmd{\ALG@doentity}{\noindent\hskip\ALG@tlm}{\ALG@printindent}{}{\errmessage{failed to patch}}
\makeatother

\algrenewcommand\algorithmicend{\strut\textbf{end}}
\algrenewcommand\algorithmicdo{\strut\textbf{do}}
\algrenewcommand\algorithmicwhile{\strut\textbf{while}}
\algrenewcommand\algorithmicfor{\strut\textbf{for}}
\algrenewcommand\algorithmicforall{\strut\textbf{for all}}
\algrenewcommand\algorithmicloop{\strut\textbf{loop}}
\algrenewcommand\algorithmicrepeat{\strut\textbf{repeat}}
\algrenewcommand\algorithmicuntil{\strut\textbf{until}}
\algrenewcommand\algorithmicprocedure{\strut\textbf{procedure}}
\algrenewcommand\algorithmicfunction{\strut\textbf{function}}
\algrenewcommand\algorithmicif{\strut\textbf{if}}
\algrenewcommand\algorithmicthen{\strut\textbf{then}}
\algrenewcommand\algorithmicelse{\strut\textbf{else}}
\algrenewcommand\algorithmicrequire{\strut\textbf{Input:}}
\algrenewcommand\algorithmicensure{\strut\textbf{Output:}}
\let\oldState\State
\renewcommand{\State}{\oldState\strut}
\algnewcommand{\IIf}[1]{\State\algorithmicif\ #1\ \algorithmicthen}
\algnewcommand{\EndIIf}{\unskip\ \algorithmicend\ \algorithmicif}

\usepackage{nomencl}
\makenomenclature
\usepackage{framed} 
\usepackage{multicol} 

\usepackage[utf8]{inputenc}
\usepackage{newunicodechar}

\DeclareUnicodeCharacter{0229}{\c{e}}
\DeclareUnicodeCharacter{2009}{ °}
\usepackage[T1]{fontenc}
\usepackage{textcomp}
\usepackage{lmodern}
\usepackage{ltxcmds}
\usepackage{gensymb}
\usepackage{siunitx}
\DeclareSIUnit{\euro}{\texteuro}
\makeatother
\definecolor{royalpurple}{rgb}{0.58, 0.44, 0.86}

\allowdisplaybreaks

\DeclareMathOperator*{\argmax}{argmax}
\DeclareMathOperator*{\argmin}{argmin}

\newfloat{subroutine}{htbp}{loa}
\floatname{subroutine}{Subroutine}

\usepackage{setspace}
\hypersetup{colorlinks = true,linkcolor = blue,anchorcolor =red,citecolor = blue,filecolor = red,urlcolor = red, pdfauthor=author}

\journal{}

\pgfplotsset{compat=1.17} 
\begin{document}

\begin{frontmatter}

\title{Modelling and analysis of multi-timescale uncertainty in energy system planning}


\author[a]{Hongyu Zhang\corref{cor1}}
\ead{hongyu.zhang@ntnu.no}
\author[a]{Erlend Heir}
\ead{asbjorn.nisi@ntnu.no}
\author[a]{Asbjørn Nisi}
\ead{erlend.heir@ntnu.no}
\author[a]{Asgeir Tomasgard}
\ead{asgeir.tomasgard@ntnu.no}

\cortext[cor1]{Corresponding author}

\address[a]{Department of Industrial Economics and Technology Management, Norwegian University of Science and Technology, Høgskoleringen 1, 7491, Trondheim, Norway}

\begin{abstract}

Recent developments in decomposition methods for multi-stage stochastic programming with block separable recourse enable the solution to large-scale stochastic programs with multi-timescale uncertainty. Multi-timescale uncertainty is important in energy system planning problems. Therefore, the proposed algorithms were applied to energy system planning problems to demonstrate their performance. However, the impact of multi-timescale uncertainty on energy system planning is not sufficiently analysed. In this paper, we address this research gap by comprehensively modelling and analysing short-term and long-term uncertainty in energy system planning. We use the REORIENT model to conduct the analysis. We also propose a parallel stabilised Benders decomposition as an alternative solution method to existing methods. We analyse the multi-timescale uncertainty regarding stability, the value of the stochastic solution, the rolling horizon value of the stochastic solutions and the planning decisions. The results show that (1) including multi-timescale uncertainty yields an increase in the value of the stochastic solutions, (2) long-term uncertainty in the right-hand side parameters affects the solution structure more than cost coefficient uncertainty, (3) parallel stabilised Benders decomposition is up to 7.5 times faster than the serial version.

\end{abstract}

\end{frontmatter}


\section{Introduction}
\label{sec:introduction}
Managing multi-timescale uncertainty is important for infrastructure planning problems. Long-term energy system planning is a type of infrastructure planning problem that is key to net-zero energy transition. Uncertainty in this type of problem comes from multiple timescales, normally including long-term and short-term timescales \citep{Kaut2014,lara2020}. Long-term uncertainty concerns timescales with years or decades, whereas short-term uncertainty relates to hourly or higher time resolution. Recently, uncertainty on the tactical timescale has been modelled and investigated for systems with significant seasonal storage \citep{Hummelen2024ExploringStorage}. In this paper, we focus on multi-timescale uncertainty that includes short-term and long-term timescales and its impact on energy system planning problems. 

The state-of-the-art modelling approach for managing multi-timescale uncertainty is Multi-Horizon Stochastic Programming (MHSP) \citep{Kaut2014}, which is essentially a type of multi-stage stochastic programming with block separable recourse \citep{Louveaux1986MultistageRecourse}. The idea of MHSP is to reduce the problem size by partially disconnecting the short-term and long-term nodes. MHSP  was used for problem with only short-term uncertainty \citep{Backe2022EMPIRE:Analyses} and problems with multi-timescale uncertainty \citep{Hellemo2013Multi-StagePerspective,Reiten2018AUncertainty, Zhang2023IntegratedDecomposition}.

The problems can still be intractable when multi-timescale uncertainty is included using an MHSP approach. To address the computational difficulty and take full advantage of MHSP, several decomposition algorithms were proposed \citep{Zhang2024DecompositionProgramming, Zhang2023IntegratedDecomposition, Zhang2024AUncertainty, Reiten2018AUncertainty}. \cite{Zhang2024DecompositionProgramming} first established that MHSP can be decomposed by Benders-type and Lagrangean-type decomposition algorithms and proposed parallel Lagrangean decomposition with primal reduction. \cite{Reiten2018AUncertainty} proposed to solve MHSP using the progressive hedging algorithm. \cite{Zhang2023IntegratedDecomposition, Zhang2024AUncertainty} proposed stabilised versions of Benders decomposition with adaptive oracles for large-scale optimisation problems with a column bounded block diagonal structure. The idea of the adaptive oracles is to avoid solving all subproblems exactly at every iteration to improve efficiency. The proposed algorithms were applied to large-scale MHSP problems with multi-timescale uncertainty and showed significant computational improvement. We extend the literature to utilise parallel computing in Benders decomposition for MHSP.

Stabilisation is very effective in addressing the oscillation issue of Benders-type decomposition when solving highly degenerate models \citep{Zhang2024AUncertainty}, which is the case for multi-region energy system planning problems. \cite{Zhang2024AUncertainty} proposed to stabilise the adaptive Benders decomposition using the level method, and a centred point stabilisation was proposed as an improvement to solve large problems with integer variables in the master problems. The centred point stabilisation was then adopted by \cite{Pecci2024RegularizedModels}
and shown significant improvement in performance. In this paper, we adopt the stabilisation method from \cite{Zhang2024AUncertainty} since our model is a large-scale multi-region model. 

Although \cite{Zhang2024AUncertainty,Zhang2023IntegratedDecomposition, Zhang2024DecompositionProgramming} proposed several efficient solution methods for MHSP, they focused on the methodology development and computational performance. The impact of multi-timescale uncertainty on energy system planning decisions is not sufficiently analysed. Therefore, in this paper, based on the methodology foundation built in \cite{Zhang2024AUncertainty, Zhang2023IntegratedDecomposition, Zhang2024DecompositionProgramming}, we rigorously analyse the impact of short-term and long-term uncertainty in energy system planning problems.

We apply the extended REORIENT model and the proposed solution algorithm to an integrated European energy system planning problem under uncertainty \citep{Heir2024NavigatingContext}. The results show that (1) the long-term uncertainty in the right hand side parameter has a significant impact on the solution structure, (2) the cost coefficient uncertainty does not change the solution structure but the magnitude, (3) the parallel stabilised Benders is up to 7.5 times faster than the serial version.

The contributions of the paper are the following: (1) we first conduct a comprehensive analysis on multi-timescale uncertainty in energy system planning, (2) we propose and test the parallel stabilised Benders decomposition, and (3) we apply the model and the algorithm to a large-scale planning problem for the European energy system to analyse the planning decisions and costs.

The outline of the paper is as follows: Section \ref{sec:multi_horizon_programming} provides the node formulation of MHSP. Section \ref{sec:Parallel stabilised Benders decomposition} proposes the parallel stabilised Benders decomposition. \ref{sec:problem_description_modelling_assumptions} presents a problem and model description. Section \ref{sec:scenario generation methods} introduces the SGRs and uncertainty assessment methods. Section \ref{sec:results} reports the computational results and numerical analysis. Section \ref{sec:conclusions} concludes the paper and suggests further research.

\section{MHSP}
\label{sec:multi_horizon_programming}
MHSP can be formulated in a node formulation or a scenario formulation, and it leads to different ways of decomposing such problems \cite{Zhang2024DecompositionProgramming}. In this paper, we focus on the node formulation and the according Benders decomposition. In the following, we present the mathematical formulation of MHSP.

\begin{equation}
   \text{MP}:\phantom{ab} \min_{\mathbf{x} \in \mathcal{X}} f(\mathbf{x})+\sum_{i \in \mathcal{I}}\pi_{i}g(\mathbf{x}_{i}, \mathbf{c}_{i}),
    \label{eq:MP}
\end{equation}
where $f(\mathbf{x})=\sum_{i \in \mathcal{I}}\pi_{i}\mathbf{c}^{\top}_{i}\mathbf{x}_{i}$ and the function $g(\mathbf{x}_{i},\mathbf{c}_{i})$ is the optimal solution of the linear programming subproblem, 
\begin{equation}
    \phantom{ab} g(\mathbf{x}_{i},\mathbf{c}_{i}):=\min_{\mathbf{y}_{i} \in \mathcal{Y}}\{\left(\mathbf{c}_{i}^{\top}C+\mathbf{c}^{\top}\right)\mathbf{y}_{i}| A\mathbf{y}_{i}\leq \mathbf{b} + B\mathbf{x}_{i}\}.
    \label{eq:SP}
\end{equation}

\begin{figure*}[!htb]
    \centering
    \includegraphics[scale=0.8]{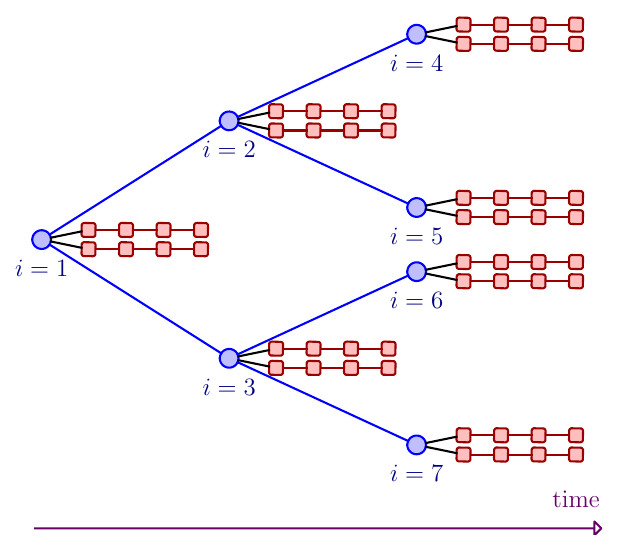}
    \caption{Illustration of MHSP with short-term and long-term uncertainty. (blue circles: strategic nodes, red squares: operational periods, $i:$ index of the strategic nodes)}
    \label{fig:MHSP}
\end{figure*}

One can have a much smaller model by disconnecting operational nodes between successive planning stages and embedding them into their respective strategic nodes. An illustration of MHSP with short-term and long-term uncertainty is shown in Figure \ref{fig:MHSP}. We call an operational problem, each block of red circles in Figure \ref{fig:MHSP}, embedded in a strategic node an operational node. The resulting model is called MHSP. MHSP is identical to multi-stage stochastic programming provided two requirements are met \citep{Kaut2014}: (a) strategic and operational uncertainties are independent, and the strategic decisions must not depend on any particular operational decisions, and (b) the operational decisions in the last operational period in a stage do not affect the system operation in the first operational period in the next stage. If either of these conditions is not met then MHSP gives only an approximation. MHSP has been applied in several energy systems planning problems \citep{Skar2016,Backe2022EMPIRE:Analyses,Zhang2022OffshoreShelf,Durakovic2023PoweringPrices}. Furthermore, the bounds in MHSP have been studied in \cite{Maggioni2020BoundsPrograms}. 

\section{Parallel stabilised Benders decomposition}
\label{sec:Parallel stabilised Benders decomposition}
Benders decomposition was first developed in \cite{Benders1962} and has been successfully applied to a wide range of difficult optimisation problems \citep{Rahmaniani2017}. Benders decomposition exploits the block diagonal structure of Equations \eqref{eq:MP}-\eqref{eq:SP} and creates outer linearisation. This method has been used in stochastic programming and is known as the L-shaped method \citep{VanSlyke1969}. Directly solving \eqref{eq:MP} -\eqref{eq:SP} can be prohibitive if there are a large number of decision nodes, which occurs in stochastic programming when there are multiple stages or many uncertain parameters.

In Benders decomposition, instead of solving Equations \eqref{eq:MP} -\eqref{eq:SP} directly, a sequence of approximations is solved. Two types of constraints can be added after each solve: feasibility cuts (enforcing the feasibility of \eqref{eq:MP}) and optimality cuts (linear approximations to \eqref{eq:MP})\citep{Birge2011}. At iteration $j$, the Relaxed Master Problem (RMP) is
\begin{subequations}
    \begin{alignat}{3}
        &\min_{\mathbf{x} \in \mathcal{X}, \beta} && f(\mathbf{x}) + \sum_{i \in \mathcal{I}}\pi_{i}\beta_{i}\\
        &\text{s.t. } && \beta_{i} \geq \theta + \boldsymbol{\lambda}^{\top}(\mathbf{x}_{i}-\mathbf{x}),&\hspace{1cm}  (\mathbf{x},\theta,\boldsymbol{\lambda}) \in \mathcal{F}_{i(j-1)}, i \in \mathcal{I},
        \end{alignat}
        \label{eq:RMP}
\end{subequations}
where $\mathcal{F}_{i(j-1)}$ is the set of cuts associated with subproblem $i$ generated prior to iteration $j$. In iteration $j$ of Benders decomposition, we first solve the RMP to obtain a solution $\mathbf{x}_j$. Then we extract the subvector $\mathbf{x}_{ij}$ of $\mathbf{x}_j$, corresponding to subproblem $i$ as its right hand side parameters. Solving this gives the optimal value of the subproblem, $\theta_{ij}$, and a subgradient, $\boldsymbol{\lambda}_{ij}$ at $\mathbf{x}_{ij}$.  Finally, new cutting planes are added to $\mathcal{F}_{i(j-1)}$ which give $\mathcal{F}_{ij}:=\mathcal{F}_{i(j-1)}\cup \{\mathbf{x}_{ij},\theta_{ij},\boldsymbol{\lambda}_{ij}\}$. This version is referred to as multi-cut Benders decomposition. The algorithm iterates until the upper bound and the lower bound converge. The lower bound is the optimal value of the RMP. The upper bound is the best feasible solution. Benders decomposition converges in a finite number of iterations when the subproblem is linear programming. In this paper, we also adopt the centred point stabilisation in \cite{Zhang2023IntegratedDecomposition} to solve large problem instances.


\begin{algorithm}[!htb]
\caption{Parallel stabilised Benders decomposition}\label{alg:stand_benders}
\begin{algorithmic}[1]
     \State choose $\epsilon$ (convergence tolerance), $\gamma$ (stabilisation factor), $\underline{\beta}$ (initial lower bound for $\beta_{i}$), $U^{*}_0:=M$ (initial upper bound), set  $j:=0$, $\mathcal{F}_{i0}:=\{(\beta_{i0},0,0)\}$ for each $i \in \mathcal{I}$;
     \Repeat
      \State $j:=j+1$;
      \State solve RMP, Equation \eqref{eq:RMP}, and obtain $\beta_{ij}$ and $\mathbf{x}^{RMP}_{j}$; set $L_{j}:=f(\mathbf{x}^{RMP}_{j})+\sum_{i \in \mathcal{I}}\pi_{i}\beta_{ij}$;
        \If{$j=1$}{}
        \StatePar{$\mathbf{x}^{Eva}_j:=\mathbf{x}^{RMP}_j$;}
        \Else{}
        \StatePar{$T_j:=L_j-\gamma\left(U^{*}_{j-1}-L_{j}\right)$;}
        \StatePar{solve centre point stabilisation problem and obtain $\mathbf{x}^{CP}_{j}$, $\mathbf{x}^{Eva}_j:=\mathbf{x}^{CP}_j$ ;}
        \EndIf
      \For{$i \in \mathcal{I}$}
      \StatePar{solve subproblem $i$, Equation \eqref{eq:SP}, using a computer node at $(\mathbf{x}^{Eva}_{ij},c_{i})$ and  obtain $\theta_{ij}$ and $\boldsymbol{\lambda}_{ij}$;}
      \EndFor
        \For{$i \in \mathcal{I}$}
      \StatePar{$\mathcal{F}_{ij}:=\mathcal{F}_{i(j-1)}\cup \{(\mathbf{x}^{RMP}_{ij},\theta_{ij},\boldsymbol{\lambda}_{ij})\}$;}
      \EndFor
    \StatePar{$U^{*}_{j}:=\min(U_{j-1},f(\mathbf{x}^{RMP}_{j})+\sum_{i \in \mathcal{I}}\pi_{i}\theta_{ij})$;}
     \Until{$U_{j}-L_{j} \leq \epsilon.$}
\end{algorithmic}
\end{algorithm}

\begin{figure}[!htb]
    \centering
    \includegraphics[scale=0.4]{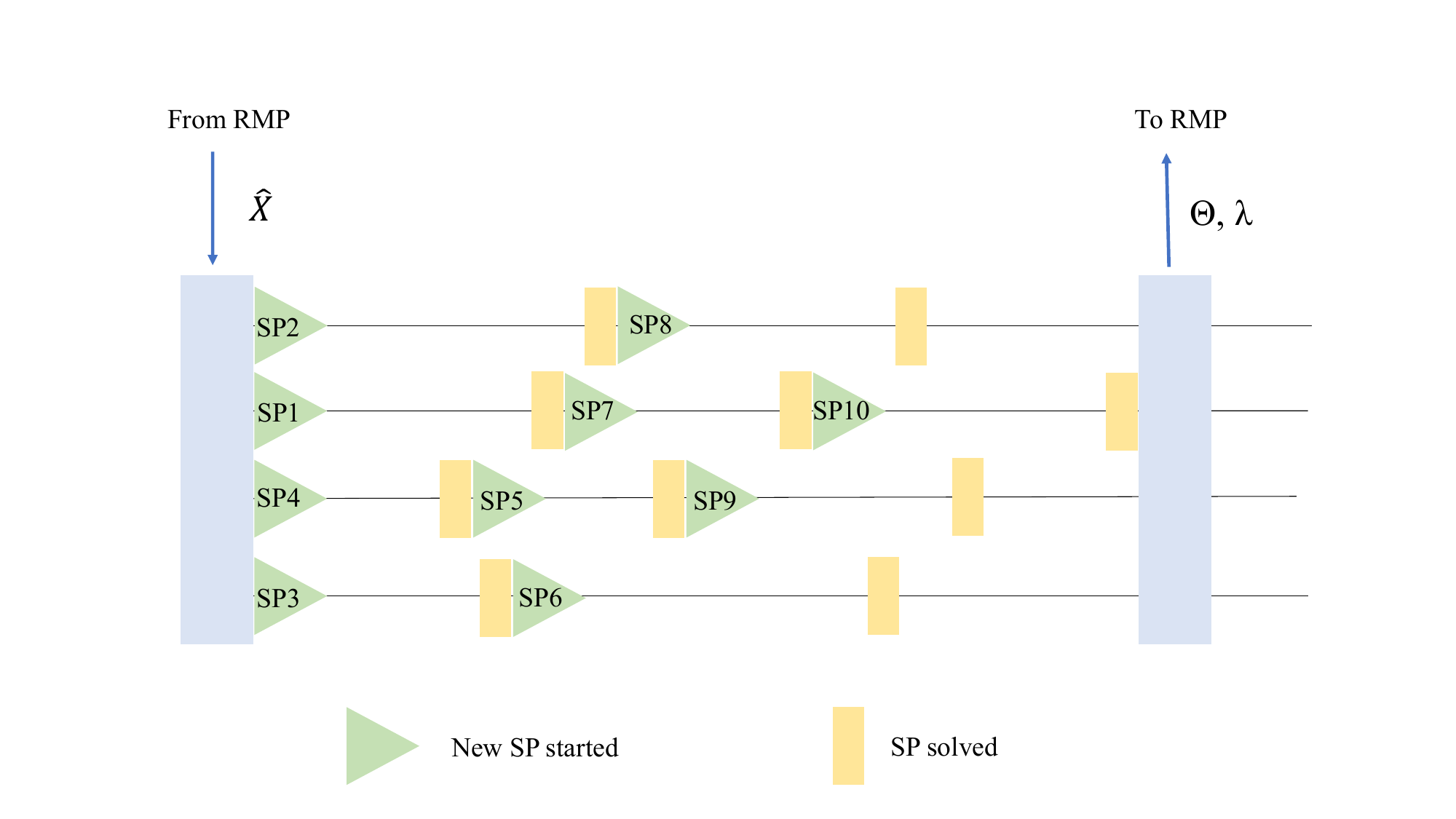}
    \caption{Whenever the computer resource is available, we start to solve a new subproblem. When all subproblems are solved, the information is fed back to the RMP.}
    \label{fig:multi-thread}
\end{figure}

The subproblems, Equation \eqref{eq:SP}, can be solved in parallel. A parallelised implementation of this algorithm is executed synchronously. The subproblems are solved in parallel whenever computational power is available, and all subproblems are solved and cuts generated before the algorithm moves on. The synchronisation of the parallellisation scheme is illustrated in Figure \ref{fig:multi-thread}. 

Parallelisation can refer to either multithreading or multiprocessing.  In multithreading, the computational workload is assigned to different logical cores on a single computer. Multiprocessing, on the other hand, assigns the workload to separate computers that communicate over a network connection. An important difference between these two methods of parallelisation is the way memory is organised. In multithreading, the computations share a common memory space on the computer they run on, while in multiprocessing, the computations have entirely separate memory spaces. This makes multiprocessing more flexible, as work can be distributed to computers that are logically and geographically separated as long as they can communicate over a network connection. Although multithreading is limited to using the resources of one computer, multiprocessing usually comes with more overhead due to network communication.

\begin{figure}[!htb]
    \centering
    \includegraphics[scale=0.4]{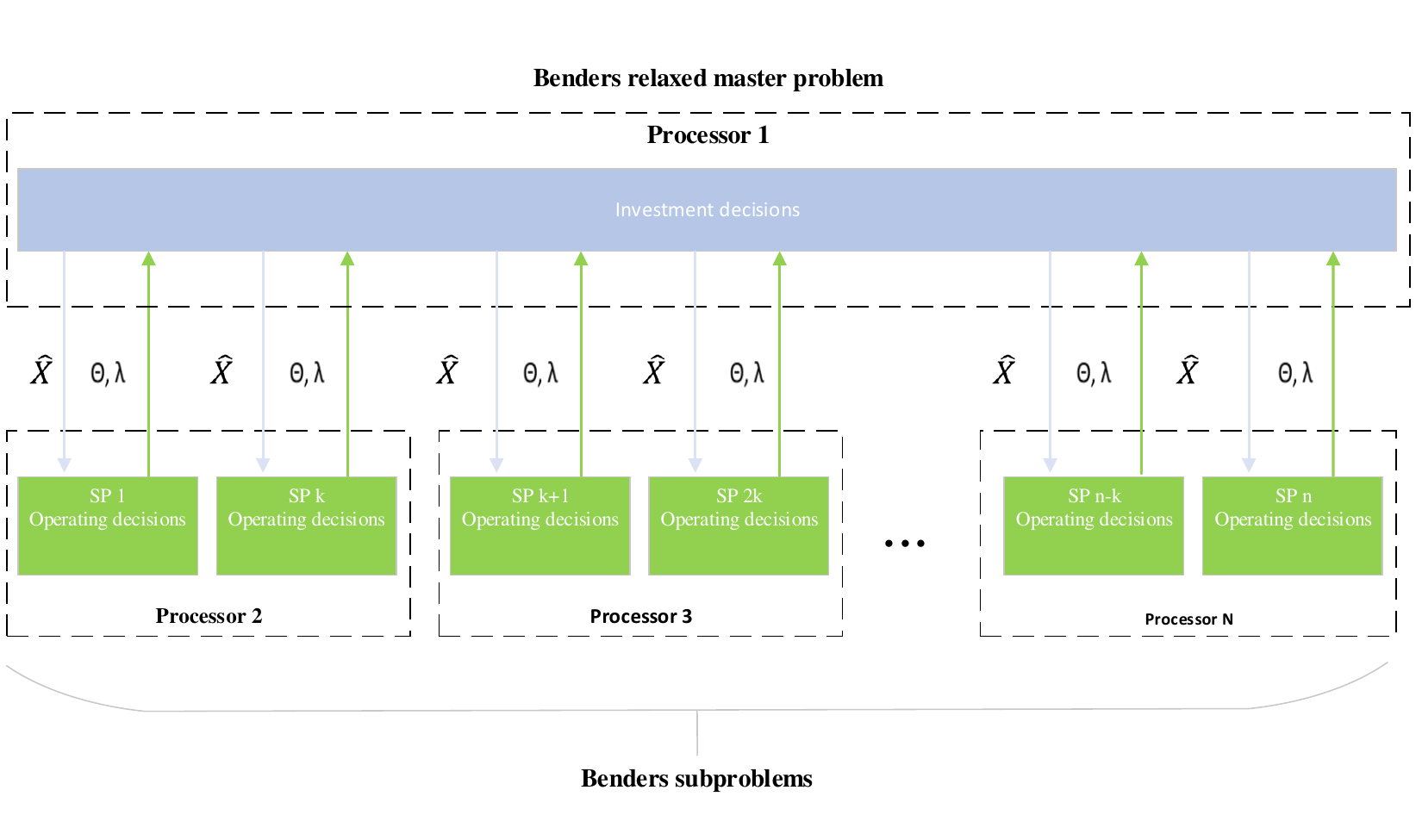}
    \caption{An overview of distributed Benders subproblems on several processors in parallel. Each process handles several subproblems and solves them on its own threads. The processors work synchronously with respect to the MP.}
    \label{fig:parallel_computing}
\end{figure}

This paper uses a parallelisation scheme that combines multithreading and multiprocessing to maximise the efficiency of the computing resources. This is orchestrated through a master processor that transfers data to multiple processors and initiates the computations. Each processor then executes the assigned computations on its own machine resources before sending results back to the master process. In our case, the subproblems are first distributed on separate computers and then solved in parallel in a multithreaded fashion. The communication between the processors is illustrated in Figure \ref{fig:parallel_computing}.

A simple static heuristic assigns computing resources to each subproblem. First, the subproblems are equally distributed on the available computers and the available logical cores on each computer are evenly distributed among the subproblems assigned to that computer. 

\section{Problem and model description}
\label{sec:problem_description_modelling_assumptions}
In this section, we present the energy system problem description. The standard assumptions and modelling strategies, and temporal and geographical representations of the problem from the original REORIENT model are kept, and we refer to \cite{Zhang2023IntegratedDecomposition} for details.

The problem under consideration aims to choose (a) the optimal strategy for investment, abandonment and retrofit planning, and (b) operating scheduling for an energy system to achieve emission targets at minimum overall costs under short-term uncertainty, including renewable energy availability, hydropower production profile and load profile, and long-term uncertainty, including oil and gas prices. 

For the investment planning, we consider: (a) thermal generators (Coal-fired plant, OCGT, CCGT, Diesel, nuclear plants, co-firing biomass with 10\% lignite, lignite); (b) generators with Carbon Capture and Storage (CCS) (Coal-fired plant with CCS and advanced CCS, gas-fired plant with CCS and advanced CCS, co-firing biomass with 10\% lignite with CCS, lignite with CCS, lignite with advanced CCS); (c) renewable generators (onshore and offshore wind and solar, wave, biomass, run-of-the-river hydropower, geothermal and regulated hydropower); (d) electric storage (hydro pump storage and lithium); (e) onshore and offshore transmission lines; (f) onshore and offshore clean energy hubs (electrolyser, fuel cell, hydrogen storage); (g) onshore steam reforming plant with CCS (SMRCCS) and (h) offshore and onshore hydrogen pipelines. The capital costs and fixed operational costs coefficients are assumed to be known. 

For the retrofit planning, we consider: (a) retrofitting existing natural gas pipelines for hydrogen transport, and (b) retrofitting existing offshore platforms for clean OEHs. Finally, we consider the abandonment of mature fields.  The problem is to determine: (a) the capacities of technologies and retrofit and abandonment decisions, and (b) operational strategies that include scheduling of generators, storage and approximate power flow among regions to meet the energy demand which minimises the combined overall expected investment and operational and environmental costs.

For a complete overview of mathematical formulation, sets and indices, parameters and variables in the REORIENT model, we refer to \cite{Zhang2023IntegratedDecomposition}. 

\section{Scenario generation and assessment methods}
\label{sec:scenario generation methods}
This paper focuses on analysing the impact of multi-timescale uncertainty in energy system planning problems. We extend the REORIENT model for this purpose. In the following, we focus on the scenario generation approaches that have been added to the REORIENT model to facilitate the analysis. 

\subsection{Short-term uncertainty scenario generation}
We adopt two scenario generation methods for short-term uncertainty: the random sampling method and the moment matching method. 

\subsubsection{Random sampling scenario generation}
The random sampling method samples scenarios from historical time series data.  To preserve both the cross- and the auto-correlation, scenarios are built using the consecutive hours from the same time slices for each data series in each region. An important factor in deciding to what degree these correlations are preserved is the length of the sequences of consecutive hours in each time slice. The random sampling method is constrained to ensure the representation of an entire year and to guarantee the capture of some seasonal variability. 

The algorithm first randomly samples a full year of data. The dataset is divided into four seasons $\mathcal{S}^{R}$, representing winter, spring, summer, and fall. The algorithm selects an equally long sequence of hours from each season $l$, ensuring that they are proportionally represented in the scenarios. These sequences represent the regular operation of the power system and are named regular seasons. In addition to regular seasons, one or more peak seasons $\mathcal{S}^{P}$ are added to each scenario. Peak seasons consist of the hours $\hat{l}$ when the grid load is at peak levels, either on aggregated or single load levels. The peak seasons aim to capture periods where the system is under larger stress and ensure that the installed capacity can handle the load levels and the ramping that may be required. The first peak season is based on the peak period for the system as a whole, while the other peak seasons are formed by using the peak hours in one area only. The areas are sorted by their peak values, and the area with the highest peak value is chosen for peak season 2, the second highest for peak season 3 and so forth. The seasons are not connected, meaning there is no storage transfer between them. 

For the later long-term periods, values are scaled from the sampled values to reflect the
predicted long-term development. The scenarios for later periods are based on the same time series as the first period. This implies that the model assumes the patterns or statistical properties in the time series are the same in the future as they are in the sampling space. This preserves the stability of the stochastic solutions. The random sampling scenario generation method is presented in Algorithm \ref{alg:Random_sampling_alg}.

\begin{algorithm}[H]
\caption{Random sampling algorithm}\label{alg:Random_sampling_alg}
\begin{algorithmic}[1]
\For {each short-term scenario $\omega \in \Omega$}
\State {Subroutines \ref{alg:subroutine} and \ref{alg:subroutine2};}
\EndFor
\end{algorithmic}
\end{algorithm}

\begin{subroutine}[!htb]
    \caption{Generate scenarios for regular seasons}\label{alg:subroutine}
    \begin{algorithmic}[1]
    \State {Select a random year $y \in \mathcal{Y}$;}
    \For {$s \in \mathcal{S}^{R}$}
    \State {Sample a random number $h$ in $[1,\frac{|\mathcal{H}|}{|\mathcal{S}^{R}|}-l-1]$;}
    \State {Set $\boldsymbol{\xi}_{s\omega}:= \boldsymbol{\xi}^{Raw}_{[h,h+l-1]sy}$;}
    \EndFor
    \end{algorithmic}
\end{subroutine}

\begin{subroutine}[!htb]
    \caption{Generate scenarios for peak seasons}\label{alg:subroutine2}
    \begin{algorithmic}[1]
        \If {$|\mathcal{S}^{P}| \geq 1$}
        \State {Set $\bar{h}:= \argmax_{h \in \mathcal{H}} \sum_{n \in \mathcal{N}}  \xi^{load}_{nhy}$;}
        \StatePar{Set: $\boldsymbol{\xi}_{(|\mathcal{S}^{R}|+1)\omega}:= \boldsymbol{\xi}^{Raw}_{[h-\frac{\hat{l}}{2}, h+\frac{\hat{l}}{2}]sy}$;}
        \If {$|\mathcal{S}^{P}| \geq 2$}
        \State{Set $j:=1$;}
        \Repeat
        \State{Set $s:=|\mathcal{S}^{R}|+j+1$;}
        \State{Set $h := \argmax_{h \in \mathcal{H}} \max_{n \in \mathcal{N}} \xi^{load}_{nhy}$;}
        \State{Set $\boldsymbol{\xi}_{s\omega}:= \boldsymbol{\xi}^{Raw}_{[h-\frac{\hat{l}}{2}, h+\frac{\hat{l}}{2}]sy}$;}
        \State{Set $\mathcal{N}:=\mathcal{N} \setminus n$;}
        \Until{$s=|\mathcal{S}^R|+|\mathcal{S}^P|$}
        \EndIf
        \EndIf
    \end{algorithmic}
\end{subroutine}

\subsubsection{Moment matching scenario generation}
The moment matching SGR aims to statistically resemble the historical data, assuming that different scenario trees with moments that match those of the true distribution have a more stable behaviour. Selecting scenarios based on matching moments is in practice a scenario reduction technique, as the scenarios with the best matching moments are selected from a larger pool of candidate scenarios.

First, a set of candidate scenario trees $\mathcal{T}$ is drawn. Let $\rho_{nst}$ be the uniform and univariate distribution for the loads in node $n$, in season $s$ and candidate scenario tree $t$. Then, the first four moments, moments, mean, variance, skewness and kurtosis, $m_{inst}$, of the candidates are calculated. The $m_{inst}$ is the $i$th moment calculated from $\rho_{nst}$ for node $n$, season $s$ and candidate scenario tree $t$. 

The quality of the candidate scenario trees is measured based on the four-moment distance \cite{Backe2021StableSampling}. The moment distance is calculated as the distance between $m_{inst}$ and $M_{ins}$, where $M_{ins}$ consists of the moments of all realisations of the electricity load in node $n$, season $s$. In addition to the scenario-independent weight, we also weigh the moments using a moment-dependent weight $\alpha$. This is based on two assumptions in \cite{Kaut2021ScenarioData}: (a) the sensitivity of the optimisation model decreases with the order of the moment, and (b) mismatch in the lower moments implies an error in the evaluation of the higher moments. We use $\alpha = \{10, 5, 1, 0.5\}$ in this paper. 

The candidate scenario with the smallest four-moment distance $\Delta M_t$ is selected and used. After selecting a scenario tree, peak seasons are added using the same method, which does not affect the moment calculations of the scenarios. The moment matching algorithm is presented in Algorithm \ref{alg:Moment_matching_alg}.

\begin{algorithm}[H]
\caption{Moment matching algorithm}\label{alg:Moment_matching_alg}
\begin{algorithmic}[1]
\For {$t \in \{1,...,\mathcal{T}\}$}
    \State {Subroutine \ref{alg:subroutine}}
\EndFor
\For {$n \in \mathcal{N}, t \in \mathcal{T}, s \in \mathcal{S}, i \in \mathcal{I}$}
\State{Calculate $m_{inst}$ from $\rho_{nst}$};
\EndFor
\State{Calculate $\Delta M_t:= \sum_{n\in \mathcal{N}}\sum_{s\in \mathcal{S}} \frac{M_{1ns}}{\sum_{n \in \mathcal{N}}M_{1ns}} \sum_{i=1}^{4} \alpha_{i}\lvert\frac{m_{inst}-M_{ins}}{M_{ins}}\rvert$;}
\State {Set $t:=\argmin \Delta M_t$;}
\State{Add peak seasons using Subroutine \ref{alg:subroutine2} to scenario $t$;}
\State Return scenario $t$.
\end{algorithmic}
\end{algorithm}

Note that the moments are calculated based on the power load values but can be generalised to other parameters. The load values are used here as they are considered the single most important factor when scaling the investments in the system. 

\subsection{Long-term uncertainty scenario generation}
Long-term scenario generation is much less studied than short-term scenario generation. Policies, technology development, and market demand are subject to a high degree of uncertainty. Long-term uncertainty scenario generation in energy systems falls into two categories: statistical methods and expert opinion methods. Statistical methods are usually based on an expected trend, either from historical trends or expert projections, and then generate scenarios based on a mathematical stochastic process. When generating scenarios based on expert opinions, data and projections from acclaimed institutions are used as a starting point for creating a scenario tree. 

This paper adopts expert opinion methods for generating scenarios for long-term, uncertain parameters. The final scenario trees produced are all unbiased with respect to an expected development trajectory for all parameters. Estimates and predictions from expert opinions regarding the long-term parameters are subject to considerable uncertainty, and there can be significant differences between experts. When generating long-term scenarios, the scenario tree reflects this by capturing the expected values and variability between future predictions.

\subsection{Assessment methods of modelling short-term uncertainty}
Here, we present the assessment metrics that are used in this paper for analysing the impact of multi-timescale uncertainty.

\subsubsection{In-sample stability of scenario generation}
Stability refers to the ability of the SGR to capture the characteristics of the underlying distribution so that it produces similar results when repeatedly generating scenarios. Tests for in-sample and out-of-sample stability are used to evaluate the scenario generation quality.

An in-sample stability test evaluates whether the stochastic program produces similar optimal objective values when different scenario trees are generated using the same scenario generation approach. If scenario trees of the same size generated using the same method produce approximately the same objective values, they are referred to as in-sample stable. The in-sample tests are first assessed separately with $m$ tests with $n$ scenarios for each SGR. The stability of the same SGR is then also assessed using a different number of scenarios. For an SGR to be in-sample stable, the following relation should hold: 
\begin{equation*}
        v^{I}_{1} \approx  v^{I}_{2} \approx ... \approx v^{I}_{m},
\end{equation*}
which means that the objective value, $v_{m}$, of each candidate sample $m$, for an in-sample stable solution, should be approximately equal to the other $|m-1|$ candidate solutions. To quantify the stability level, the mean and Standard Deviation (SD) are calculated for the in-sample values. The significance of the results depends on the number of tests $m$.

\subsubsection{Out-of-sample stability of scenario generation}
Out-of-sample stability refers to whether different solutions from the same SGR give the same solution when tested on the true distribution. We first obtain a solution $\mathbf{x}^{n^\prime}$ by solving a problem with $n^\prime$ scenarios, where $n^\prime$ is as large as possible to approximate the true distribution. Then we use the same SGR to generate $m$ instances with $n$ scenarios, and we fix the first stage solution to $\mathbf{x}^{n^\prime}$ and solve the underlying instance to obtain objective values. In the context of the REORIENT model, the first-stage solutions are the investment decisions resulting from solving the model with a particular short-term scenario tree.  Out-of-sample stability is characterised by the out-of-sample values being approximately equal: 
\begin{equation*}
    v^{O}_{1}(\mathbf{x}^{n^\prime}) \approx v^{O}_{2}(\mathbf{x}^{n^\prime}) \approx ... \approx v^{O}_{m}(\mathbf{x}^{n^\prime}).
\end{equation*}

In addition to the stability of the objective function, the solution stability within and across the different methods is assessed. By solution stability, we refer to the investment decisions taken by the model. If an SGR gives the same or similar investments for each replication, it is solution stable. This is to see whether different samples generated with the same SGR produce similar solutions. In addition, there can be many degenerate solutions with equal or close to equal objective values. This is especially relevant in the REORIENT model, using a decomposition algorithm until it reaches an acceptable convergence level, thus not guaranteeing the true optimal solution. Studying how the solution varies is therefore important in order to use the model as a planning tool.

\subsubsection{Value of the Stochastic Solution (VSS) of short-term uncertainty}
We also evaluate the VSS for short-term uncertainty. The VSS is also known as the expected cost of disregarding uncertainty. The Multi-Horizon Expected Value problem (MHEV), is the MHSP equivalent of the EV. Following the methodology by \cite{Maggioni2020BoundsPrograms}, the MHEV is found by setting the strategic and operational parameters to their expected values.

The value of the stochastic solution with only operational uncertainty is calculated by first evaluating MHEV. The strategic decisions are fixed to their optimal values from the MHEV, and the model is evaluated in a setting including operational uncertainty. This solution is compared with the solution to the MHSP with operational uncertainty only. This metric quantifies the value of including stochastic short-term parameters compared with the deterministic approach. 

\subsection{Assessment methods of modelling long-term uncertainty}
To evaluate the value of including long-term uncertainty in the model, VSS and Rolling Horizon Value of Stochastic Solution (RHVSS) are used \citep{Maggioni2014BoundsProgramming}.

\subsubsection{VSS of long-term uncertainty}
To calculate the VSS of long-term uncertainty, the Expectation of the Expected Value Problem (EEV) is first calculated. This is done by first setting the stochastic parameters to their expected values, solving the deterministic problem. Then, the strategic variables are fixed to their values from the deterministic solution and evaluated on the MHSP problem. The objective value from this calculation is the EEV. The VSS is then calculated as the difference between the EEV and the objective value of the Stochastic Problem (SP). We define VSS in Equation \eqref{VSS},
\begin{align}\label{VSS}
    VSS = EEV - SP.
\end{align}

\subsubsection{RHVSS}
\begin{figure}[!htb]
    \centering
    \includegraphics[width=\textwidth]{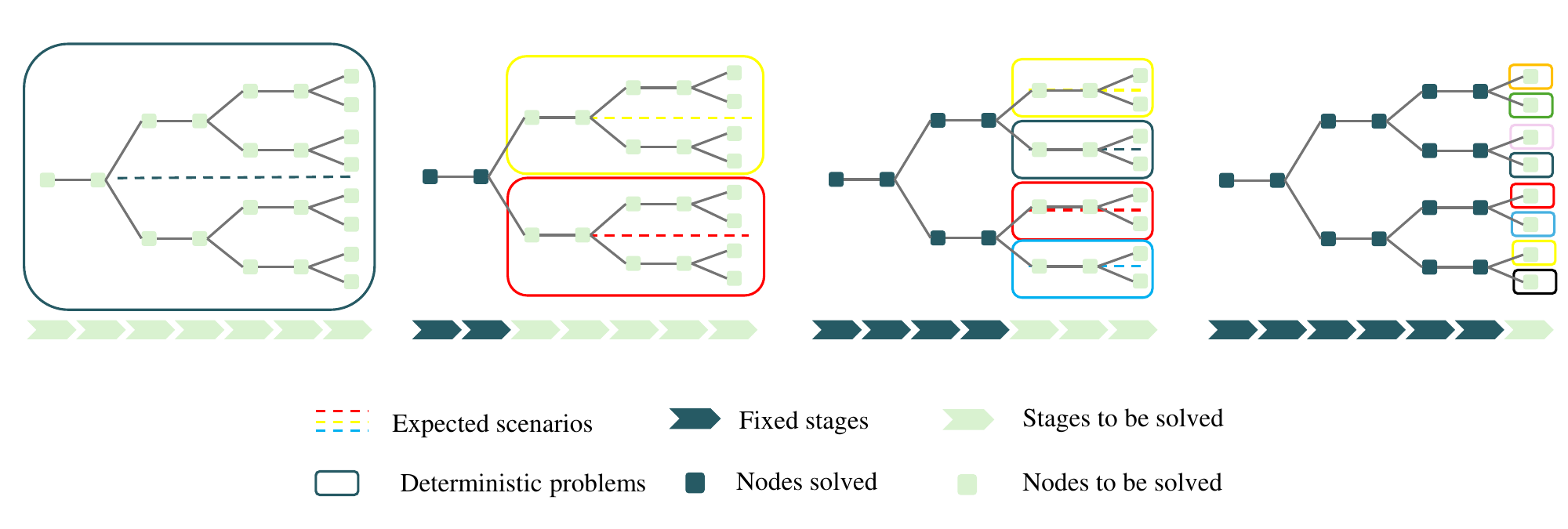}
    \caption{Illustration of the procedure for calculating the ERHEV \citep{Heir2024NavigatingContext}. The problem is solved stage by stage, fixing the solved nodes and updating the expected scenarios for the deterministic problems.}
    \label{fig:RH_scen_tree}
\end{figure}

RHVSS calculates the difference between stochastic programming and a Rolling Horizon (RH) problem. The RH problem solves a series of deterministic problems where the state is updated to reflect realised future outcomes of uncertainty. The resulting value is called the Expected Result of Using the Rolling Horizon Expected Value solution (ERHEV). Obtaining this value involves solving a deterministic problem for each branching point of the scenario tree, illustrated in Figure \ref{fig:RH_scen_tree}. 

First, the problem is solved by setting the long-term parameters to their expected values. After solving the model, the decisions made before the first stochastic branching point, in the case of Figure \ref{fig:RH_scen_tree}, the stage 1 and 2 decisions are fixed. As new information arrives between the second and third stages, this is accounted for by solving one deterministic problem for each branch of the scenario tree. These two problems share stage 1 and stage 2 decisions, but the parameter values are updated to reflect the expected value of the possible futures from the branch. The process of fixing decisions made before the next branching point prior ot solving new deterministic problems with updated parameters is repeated until one deterministic problem is solved for each scenario. The ERHEV is then obtained by calculating the weighted solutions from these solutions. 

Finally, the RHVSS can be calculated: 
\begin{equation*}
    RHVSS = MHSP - ERHEV \geq 0.
\end{equation*}

The RHVSS quantifies the value of explicitly encoding stochastic parameters in the model formulation instead of using the expected value of future outcomes. This measure reflects the fact that decision-makers are likely to update their beliefs about the future as uncertain outcomes are realised.

\section{Results}
\label{sec:results}
In this section, we first present the case study. Then we report the computational performance of the parallel stabilised Benders decomposition, followed by the analysis of the impact of multi-timescale uncertainty. 

\subsection{Case study}
\begin{figure}[!thb]
    \centering
    \includegraphics[scale=0.9]{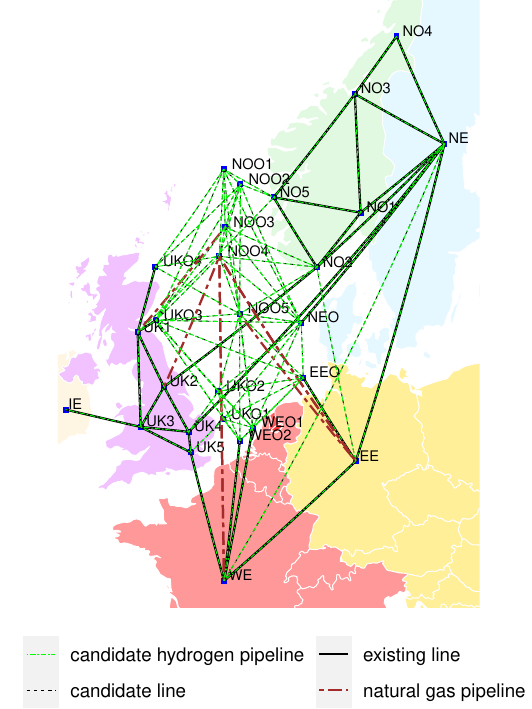}
    \caption{Illustration of the considered European energy system. The considered system includes 27 regions (each region can deploy 36 technologies), 87 transmission lines, 7 existing natural gas pipelines that can be retrofitted for hydrogen transport (some are overlapped), and 87 candidate new hydrogen pipelines.}
    \label{fig:ns_grid}
\end{figure}

We apply the extended REORIENT model to the integrated strategic planning of the European energy system \citep{Heir2024NavigatingContext}. The network topology is shown in Figure \ref{fig:ns_grid}.
We make investment planning towards 2050 with a 5-year planning step. We implemented the algorithm and model in Julia 1.8.2 \citep{Bezanson2017Julia:Computing} using JuMP \citep{jump} and solved with Gurobi 10.0 \citep{gurobi}. The problem instances contain up to 55 million continuous variables, 1876 binary variables, and 116 million constraints. We run the code on nodes of a computer cluster with a 2x 3.5GHz 8 core Intel Xeon Gold 6144 CPU and 384 GB of RAM, running on CentOS Linux 7.9.2009. 

\subsection{Computational results}
In this section, we present the computational results regarding the algorithm performance. We compare the parallel and serial implementations of the stabilised Benders decomposition: (1) a standard serial implementation, (2) a parallel multi-threaded implementation using only one computer, and (3) a parallel distributed implementation using three computers and combined with multi-threaded execution on the computer cluster.  A summary of the cases used in the computational study is presented in Table \ref{table:Overview of the cases used in the computational study.}.

\begin{table}[!htb]
    \caption{Overview of the cases used in the computational study.}
    \label{table:Overview of the cases used in the computational study.}
    \resizebox{\columnwidth}{!}{%
        \begin{tabular}
            {cS[table-format=3.0]S[table-format=1.0]S[table-format=2.0]|S[table-format=2.0]S[table-format=2.0]|S[table-format=1.1e1]S[table-format=4.0]S[table-format=1.1e1]}
            \toprule
            & {Operational periods} & \multicolumn{1}{c}{Short-term} & \multicolumn{1}{c}{Long-term} & \multicolumn{2}{c}{Number of decision nodes} & \multicolumn{3}{c}{Problem size (undecomposed)}\\
            & {per short-term scenario} & \multicolumn{1}{c}{scenarios} & \multicolumn{1}{c}{scenarios} & \multicolumn{1}{c}{Operational nodes} & \multicolumn{1}{c}{Investment nodes} & \multicolumn{1}{c}{Continuous variables} & \multicolumn{1}{c}{Binary variables} & \multicolumn{1}{c}{Constraints} \\ \hline
            Case	1	&	409	&	4	&	1	&	6	&	6	&	6.6e6	&	804	&	1.4e7	\\
            Case	2	&	985	&	4	&	1	&	6	&	6	&	1.6e6	&	804 &	3.3e7	\\
            Case	3	&	409	&	4	&	8	&	21	&   14 	&	2.3e7	&	1876	&	4.9e7	\\
            Case	4	&	985	&	4	&	8	&	21	&	14	&	5.5e7	&	1876	&	1.2e8	\\
            \bottomrule
        \end{tabular}
    }
\end{table}

From the results in Table \ref{table:Results of computational time.}, it is clear that both the multi-threaded and the distributed implementations outperform the serial implementation. Compared with the serial implementation, distributed implementation is up to 7.5 times faster, and the multi-threaded implementation is up to 3.9 times faster. The results show that an increased problem scale makes a parallel implementation more attractive, and a distributed implementation scales better than the multi-threaded implementation as the number of subproblems increases. Two notable observations are that the improvement in distributed implementation is only slightly better when compared to the available computing resources. It is also clear that multi-threading the solving of the subproblems, even on a single machine, drastically improves solving time.

\begin{table}[!htb]
    \caption{Results of computational time.}
    \label{table:Results of computational time.}
    \resizebox{\columnwidth}{!}{%
        \begin{tabular}
            {cS[table-format=4.0]S[table-format=3.0]S[table-format=4.0]|S[table-format=5.0]S[table-format=3.0]S[table-format=5.0]|S[table-format=5.0]S[table-format=3.0]S[table-format=5.0]}
            \toprule
            &  \multicolumn{3}{c}{Distributed} & \multicolumn{3}{c}{Multi-threaded} & \multicolumn{3}{c}{Serial}\\
            & {Solving time (s)} & \multicolumn{1}{c}{Overhead time (s)} & \multicolumn{1}{c}{Total (s)} & \multicolumn{1}{c}{Solving time (s)} & \multicolumn{1}{c}{Overhead time (s)} & \multicolumn{1}{c}{Total (s)} & \multicolumn{1}{c}{solving time (s)} & \multicolumn{1}{c}{Overhead time (s)} & \multicolumn{1}{c}{Total (s)} \\ \hline
            Case	1	&1826		&324		&	2150	&2373		&195		&	2569	&6722		&179 		&6901\\
            Case	2	&4805		&588		&	5393	&6254		&288		&	6543	&18781		&213 		& 18994\\
            Case	3	&2630		&357		&	2987	&5377		&360   &	5737	&22052		&209 		& 22261\\
            Case	4	&8241		&720		&	8961	&14801		&733		&	14801	&55028		&240 		& 55269\\
            \bottomrule
        \end{tabular}
    }
\end{table}

Table \ref{tab:comp_study_results_2} shows the results where we increase the number of subproblems to 105. It is clear that an increased number of computers improves solving time. In addition, a comparison between distributing work on 3 and 6 processors is performed to examine how the distributed implementation scales when more computers are used. It is clear that this improved performance does not scale linearly with the amount of computing resources.

\begin{table}[!htb]
    \caption[Computational Study - Results Machine scaling]{How the implemented algorithm scales with increasing computational resources.}
    \label{tab:comp_study_results_2}
    \centering
    \small
    \begin{tabular}{l|S[table-format=3.0]|S[table-format=5.0]|S[table-format=3.0]|S[table-format=5.0]}
        \toprule
        {\# Workers} & {Number of subproblems} & {Solving time (s)} & {Overhead time (s)} & {Total time (s)} \\
        \hline
         3 processors  &  105& 11676& 603& 12279\\
         6 processors &  105& 8790& 714 & 9504\\
        \bottomrule
    \end{tabular}
\end{table}


\subsection{Short-term uncertainty analysis}
\label{sec:Short-term uncertainty analysis}
Here, we focus on objective and solution stability analysis regarding short-term uncertainty and the evaluation of VSS. 

\subsubsection{Objective stability}
\begin{figure}[!htb]
    \centering
    \includegraphics[width=0.5\textwidth]{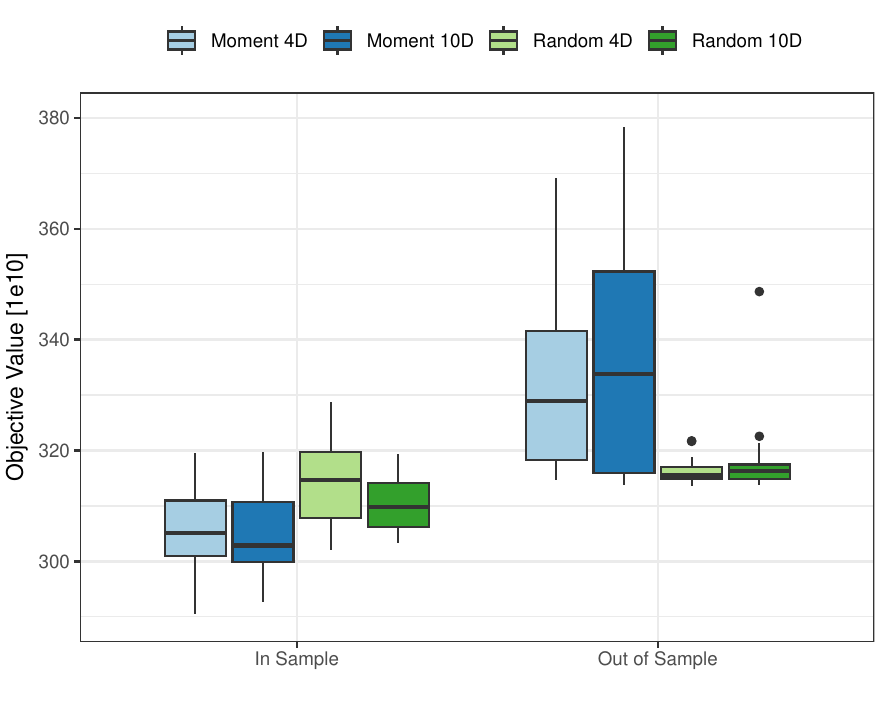}
     \caption{Comparison of stability test results for random sampling and moment matching SGRs. The boxes represent the 25\% - 75\% percentile, while the whiskers extend to the minimum and maximum values within 1.5 times the interquartile range. Outliers are plotted as individual points beyond the whiskers, and extremely large outliers are excluded from the plot. The median is represented by the line within each box.}
    \label{fig:box_in_out_of_sample}
\end{figure}

The values from the in-sample and out-of-sample stability tests of the two SGRs are presented in Figure \ref{fig:box_in_out_of_sample}. As can be observed, the random sampling and moment matching methods perform similar in-sample. Table \ref{tab:test_results_in_sample_random_7} presents the mean and SD for the in-sample tests. The different methods are observed to be relatively in-sample stable. The SD for the methods varies between 1.49\% for 10-day random to 2.45\% for 10-day moment. 

\begin{table}[!htb]
    \centering
        \caption{In sample stability.}
        \resizebox{\columnwidth}{!}{%
    \begin{tabular}{cS[table-format=2.0]S[table-format=3.0]S[table-format=1.0]S[table-format=1.0]S[table-format=2.0]S[table-format=1.2e1]S[table-format=1.2e1]r}
        \toprule
        SGR &{Num tests}& {Ses. len (h)} &{Num ses.} & {Num Peaks} & {Peak len. (h)} &  {Obj. Val. (€)} & {SD (€)} & Relative SD\\
        \midrule
        Random & 30 & 96 & 4&1 & 25 & 3.14e12 & 7.29e10& 2.32 \%\\
        Moment & 30 & 96 & 4&1 & 25 & 3.06e12 & 7.03e10  &2.30 \% \\
        Random & 30 & 240 & 4&1 & 25 & 3.10e12 & 4.61e10 &1.49 \% \\
        Moment & 30 & 240 & 4&1 & 25 &  3.03e12 & 7.47e10 &2.45 \% \\
        \bottomrule
    \end{tabular}
    }
     \label{tab:test_results_in_sample_random_7}
\end{table}

The random sampling performs significantly better than moment matching for out-of-sample stability.  The mean objective values and SD, are shown in Table \ref{tab:test_results_out_of_sample_random_7}. Moment matching with 4 days is the most unstable method out-of-sample, performing significantly worse than the other methods with an SD of 22.14\%. The 4-day random sampling performs best, having a low SD of 0.59\%. 

\begin{table}[!htb]
    \centering
        \caption{Out-of-sample stability. }
    \resizebox{\columnwidth}{!}{%
    \begin{tabular}{cS[table-format=2.0]S[table-format=3.0]S[table-format=1.0]S[table-format=1.0]S[table-format=2.0]S[table-format=1.2e1]S[table-format=1.2e1]r}
        \toprule
        SGR &{Num tests}& {Ses. len (h) } &{Num ses.} & {Num Peaks} & {Peak len (h)} &  {Obj. Val. (€)} & {SD (€)} & Relative SD\\
        \midrule
        Random & 30 & 96 & 4&1 & 25 & 3.16e12  & 1.89e10 & 0.59\%\\
        Moment & 30 & 96 & 4&1 & 25 &  3.65e12 & 8.09e11 &22.14\%\\
        Random & 30 & 240 & 4&1 & 25 & 3.18e12 & 6.15e10& 1.94\%\\
        Moment & 30 & 240 & 4&1 & 25 & 3.52e12  & 3.86e11 & 10.98\%\\
        
        \bottomrule
    \end{tabular}
    }
     \label{tab:test_results_out_of_sample_random_7}
\end{table}

\begin{figure}[!htb]
    \centering
    \includegraphics[width=0.5\textwidth]{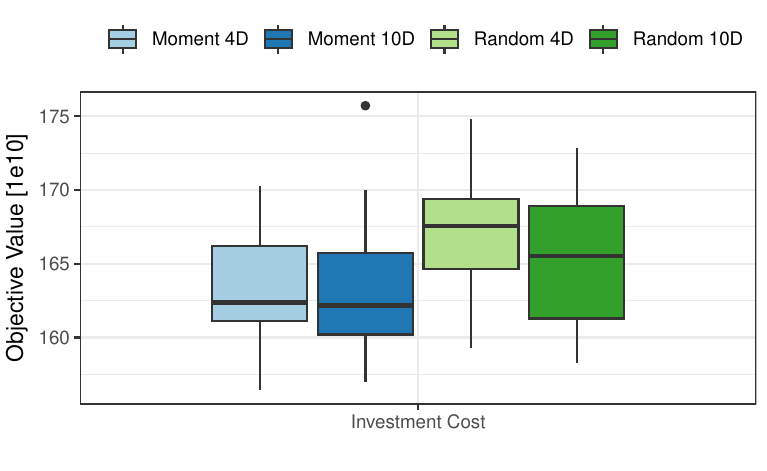}
    \caption{Comparison of investment cost for random sampling and moment matching scenario generation methods with a season length of 7 days and peak season of 25 hours.}
    \label{fig:out_of_sample_investmetn}
\end{figure}

\begin{figure}[!htb]
    \centering
    \includegraphics[width=0.5\textwidth]{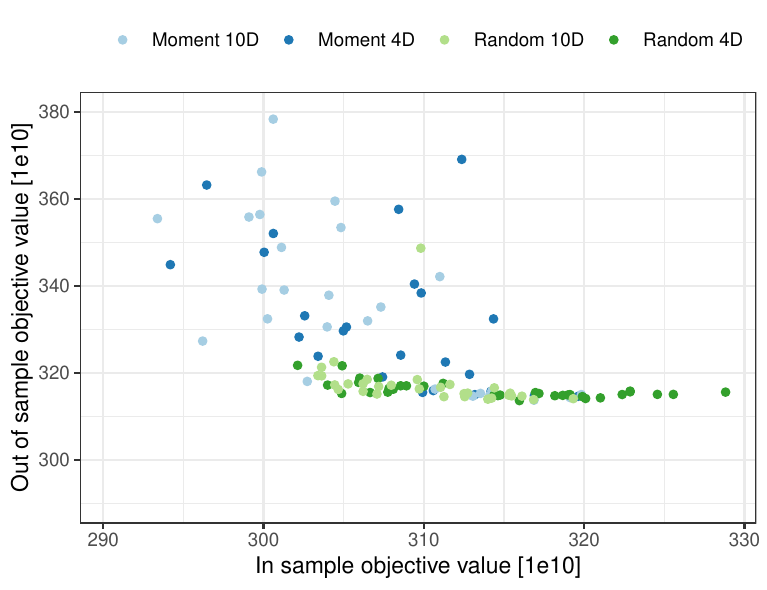}
    \caption{Comparison between the in-sample and out-of-sample values for the different SGRs. (extremely large outliers are excluded from the plot)}
    \label{fig:scatter_withput_outliers_sample_without_outliers}
\end{figure}

The objective values for the investment costs are presented in Figure \ref{fig:out_of_sample_investmetn}. As the investment decisions are fixed in the out-of-sample tests, these values are kept the same in both the in- and out-of-sample tests. The pattern here is similar to what is for the in-sample-values in Figure \ref{fig:scatter_withput_outliers_sample_without_outliers}, with the mean of the moment matching method being lower than for random sampling.  

In Figure \ref{fig:scatter_withput_outliers_sample_without_outliers}, the in-sample and out-of-sample values for the different SGRs are plotted against each other. We can see that the plot showing a clear trend that the outliers among the out-of-sample objectives tend to have a low in-sample objective value.

The moment matching values, having the lowest mean values, are also observed to be more prominently featured among the outliers and high out-of-sample values. Note that the moment matching algorithm selects the scenarios based on candidates created by random sampling. Hence, they are selected from the same distribution as the random scenarios. This observation shows the moment matching algorithm selecting a skewed selection from the random distribution.


\subsubsection{Technology stability}
\begin{figure}[!htb]
    \centering
    \includegraphics[width=0.5\textwidth]{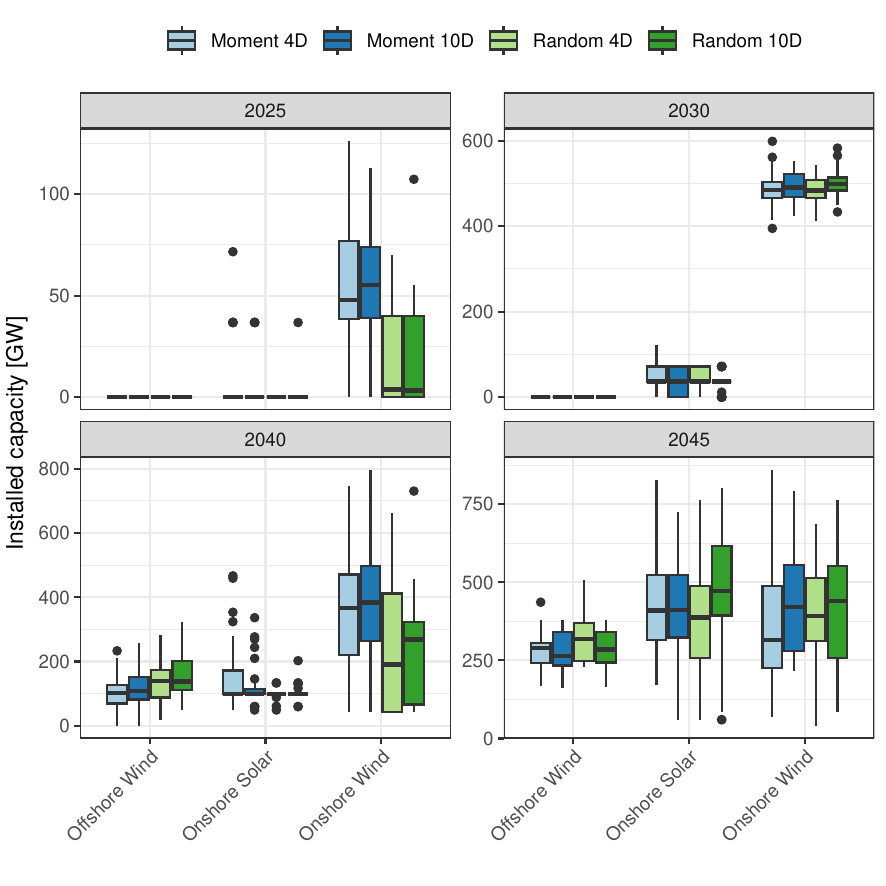}
    \caption{Comparison of period investments in renewable power generation for different SGRs (in-sample tests).}
    \label{fig:tech_stability_renewables}
\end{figure}

In addition to the objective value stability, the solution structure of the different SGRs is compared. In Figure \ref{fig:tech_stability_renewables}, the investments in the key renewable technologies are compared for the two methods. As can be observed, the variations in the technology mix are large for both methods, especially in later periods. The main outliers in the plot are in the moment matching scenarios. 

\begin{figure}[!htb]
    \centering
    \includegraphics[width=0.5\textwidth]{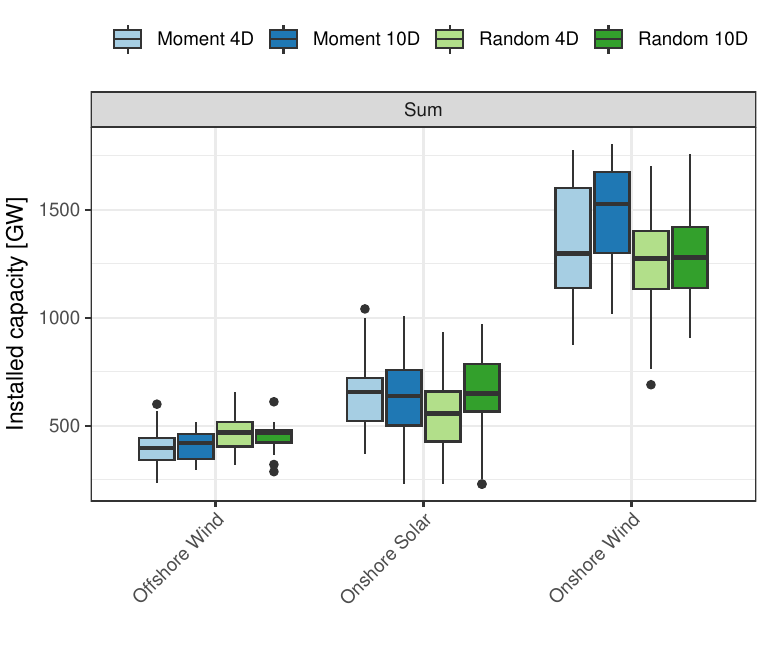}
    \caption{Comparison of the aggregated investments in renewable power generation for different SGRs.}
    \label{fig:tech_stability_renewables_aggragated}
\end{figure}

The aggregated investments into renewables for the period are presented in Figure \ref{fig:tech_stability_renewables_aggragated}. The trends in the two methods are similar, with the largest investments in onshore wind. Despite some variations in the mean values, the overall structure of the renewable investments is similar across all SGRs. This shows that systems of relatively similar composition can perform very differently with regard to stability. 

\subsection{VSS of short-term uncertainty}\label{bounds_evaluation}
The results from the evaluation of bounds for the operational uncertainty are presented in this section. The VSS is reported for different SGR configurations in Table \ref{tab:test_results_VSS_short_term}. We can see that the VSS values are large for all instances. This indicates that the first-stage solutions derived from the EV approach are not able to meet the demand in all hours. The high VSS values suggest that there are significant discrepancies between the stochastic and deterministic solutions, emphasising that the EV-based solutions are insufficient for ensuring demand fulfilment across the entire period.

\begin{table}[!htb]
    \centering
        \caption{Test Results for VSS of short-term uncertainty.}
    \resizebox{\columnwidth}{!}{%
    \begin{tabular}{cS[table-format=3.0]S[table-format=1.0]S[table-format=1.0]S[table-format=2.0]S[table-format=1.0]S[table-format=3.2]S[table-format=4.2]S[table-format=3.2]S[table-format=4.2] r}
        \toprule
        {SGR} & {Ses. len. (h)} & {Num ses.} &{Num peak.} &  {Peak len. (h)} & {Num nodes} & {EV (€bn)} & {EEV (€bn)} & {SP (€bn)} & {VSS (€bn)} & {VSS in \%}\\
        \midrule
        Random & 96 & 4 &1 & 25 & 4  & 258.36 & 2250.67 & 280.66 & 1970.01 & 87.5\%\\
        Random & 240 & 4&1 & 25 & 4 & 260.65& 1508.42& 278.03 & 1230.39 & 81.5\%\\
        Moment & 96 & 4&1 & 25 & 4  & 250.88& 906.32& 263.05& 643.27 & 71.0\%\\
        Moment & 240 & 4&1 & 25 & 4 & 265.75 &292.83 & 267.56& 25.27 & 8.6\%\\
        \bottomrule
    \end{tabular}
    }
    \label{tab:test_results_VSS_short_term}
\end{table}

It is important to note that the load shedding cost is set to a very high value in this model. This high penalty for unmet demand results in a substantial increase in the total cost when the EV solutions fail to meet demand, thereby amplifying the VSS values. 

\subsection{Long-term uncertainty analysis}
\label{sec:Long-term uncertainty analysis}
For the long-term uncertainty, we report the value of the stochastic solution via VSS and RHVSS. Also, we present the impact on the energy system solutions based on different uncertainty parameters. 

\subsubsection{Uncertainty on one parameter}\label{one_dimensional_uncertainty}
\begin{figure}[!htb]
    \centering
    \begin{subfigure}[b]{0.48\columnwidth}
        \includegraphics[width=\linewidth]{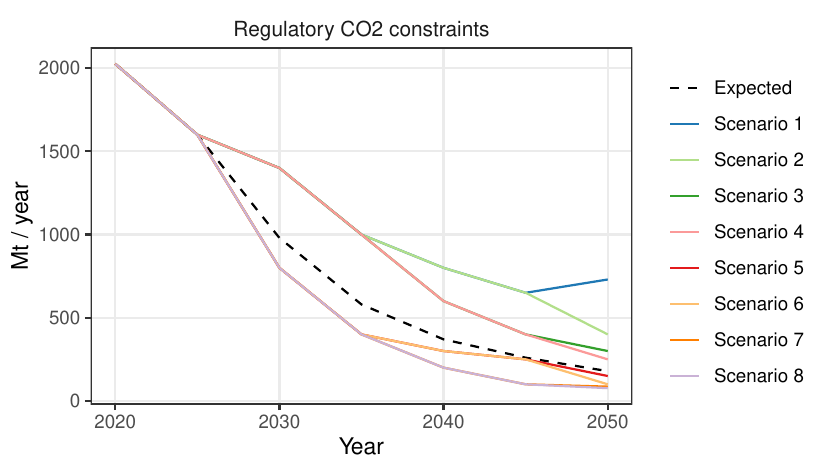}
    \end{subfigure}%
    \begin{subfigure}[b]{0.48\columnwidth}
        \includegraphics[width=\linewidth]{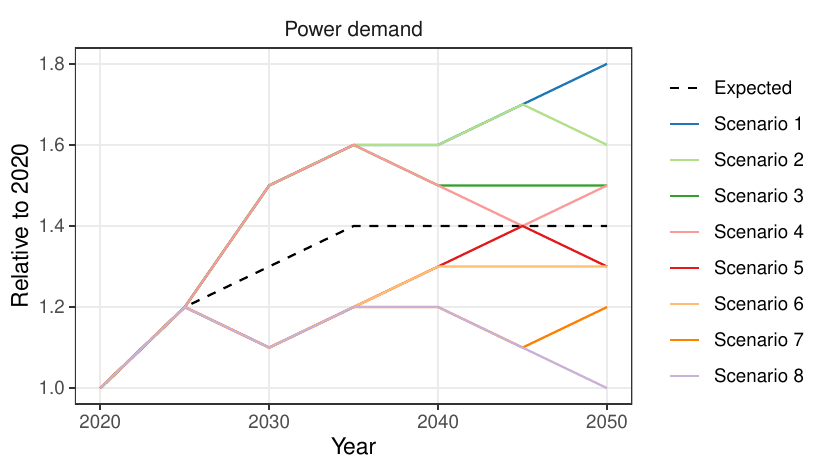}
    \end{subfigure}

    \begin{subfigure}[b]{0.48\columnwidth}
        \includegraphics[width=\linewidth]{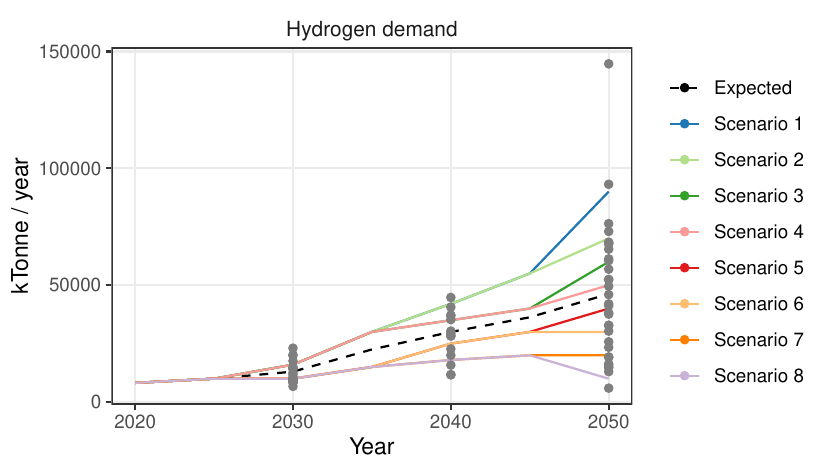}
    \end{subfigure}%
    \begin{subfigure}[b]{0.48\columnwidth}
        \includegraphics[width=\linewidth]{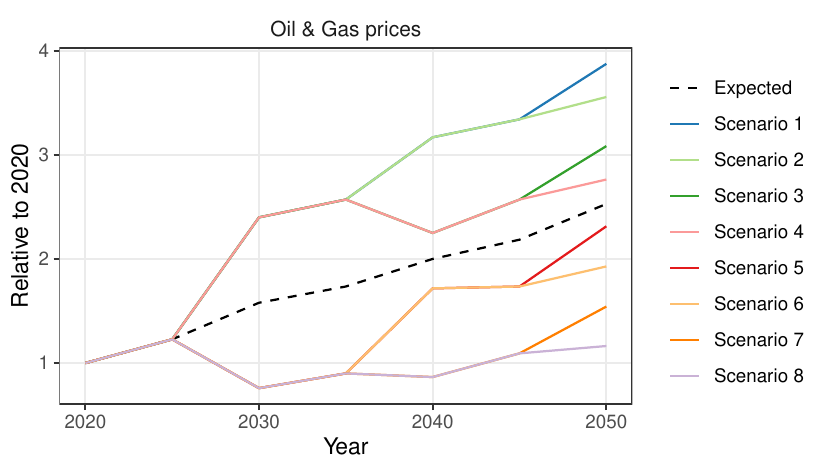}
    \end{subfigure}

    \begin{subfigure}[b]{0.48\columnwidth}
        \includegraphics[width=\linewidth]{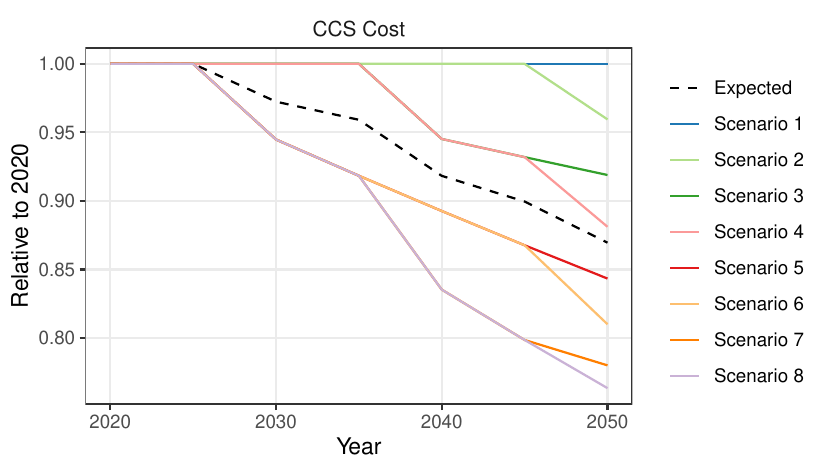}
    \end{subfigure}%
    \begin{subfigure}[b]{0.48\columnwidth}
        \includegraphics[width=\linewidth]{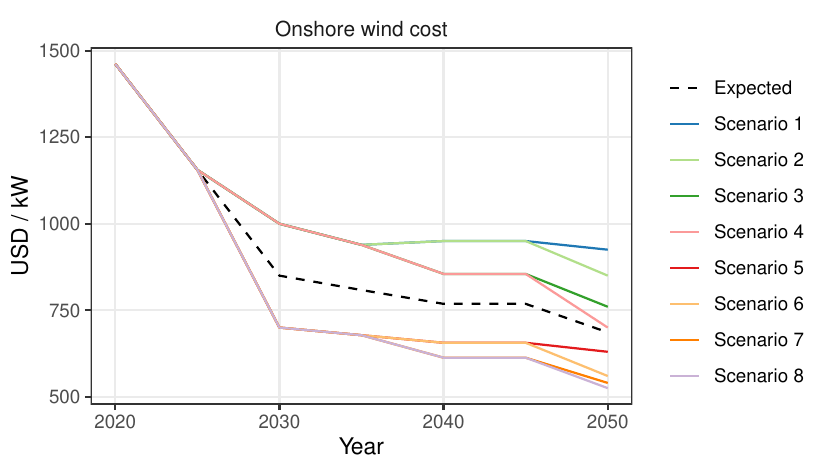}
    \end{subfigure}

    \begin{subfigure}[b]{0.48\columnwidth}
        \includegraphics[width=\linewidth]{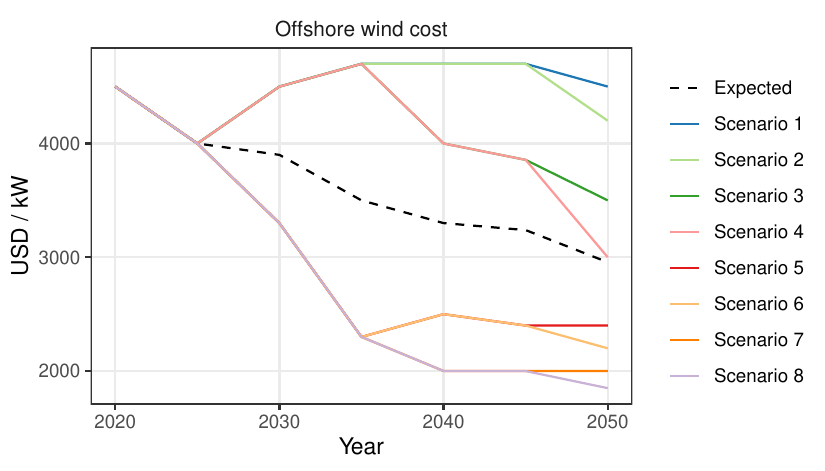}
    \end{subfigure}%
    \begin{subfigure}[b]{0.48\columnwidth}
        \includegraphics[width=\linewidth]{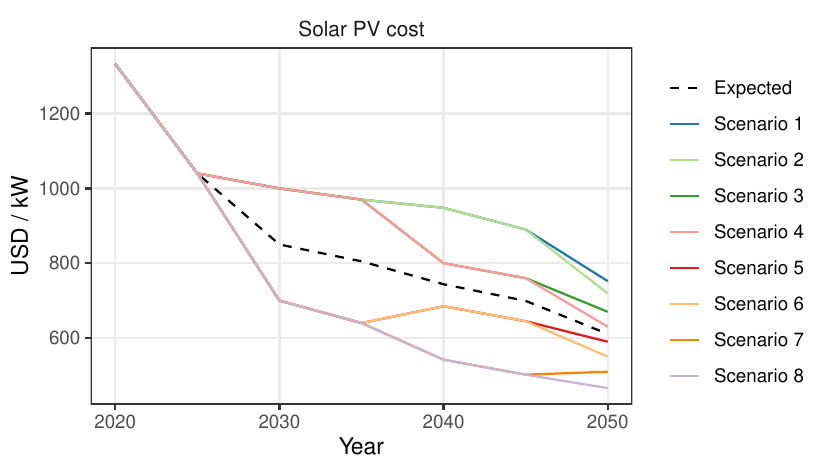}
    \end{subfigure}

    \caption{Scenario trees for each single uncertain parameter.}

    \label{fig:uncertainty scenarios}
\end{figure}

In this section, we analyse the effect of incorporating long-term uncertainty one long-term parameter at a time. We consider uncertainty in CO$_2$ cap, power demand, hydrogen demand, CCS technology costs, renewable technology costs (onshore wind, offshore wind and solar PV), and oil and gas prices. The scenario trees for these parameters are presented in Figure \ref{fig:uncertainty scenarios}.

The differences across the different parameters are first assessed by analysing the expected solutions from the MHSP in comparison with the deterministic solution. We then analyse individual scenarios. Note that the expected solutions of the MHSP are the weighted sums of the different scenarios and can be quite different from the solutions of the individual scenarios. This can give an indication of important differences between solving deterministic and stochastic models. The deterministic solution serves as a reference for comparing the investments across different solutions.

\begin{figure}[!htb]
    \centering
    \includegraphics[height=7cm]{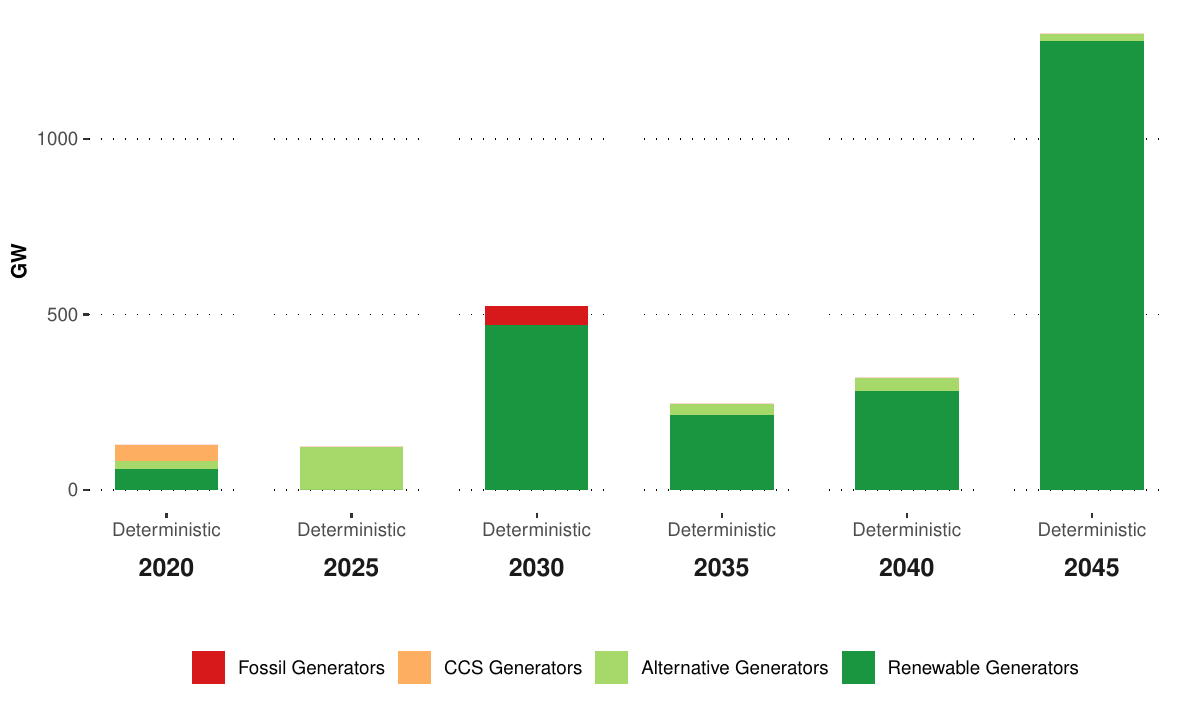}
    \caption{Investment decisions when solving REORIENT deterministically.}
    \label{fig:deterministic_generators}
\end{figure}

Figure \ref{fig:deterministic_generators} presents the investments in power generation in the deterministic case in the period 2020 to 2045. The investment decisions are aggregated into different groups of technologies, including fossil generators, CCS generators, renewables, and alternative generators. Alternative generators include nuclear, fuel cells, bioenergy, and waste-to-energy. 

The majority of investments are made in the later periods. This is because the CO$_2$ cap is reduced while demand increases simultaneously. Most new investments are in renewable generation, addressing both the CO$_2$ cap and the rising demand. Early investment periods are also characterised by investments in CCS generators and alternative generators. The investments in 2030 are significantly higher than in the preceding and following stages. This may be because the expected power demand flattens out from 2035 onwards, making it necessary to make investments to satisfy this stable but high demand level.

\begin{figure}[!htb]
    \centering
    \includegraphics[width=\textwidth]{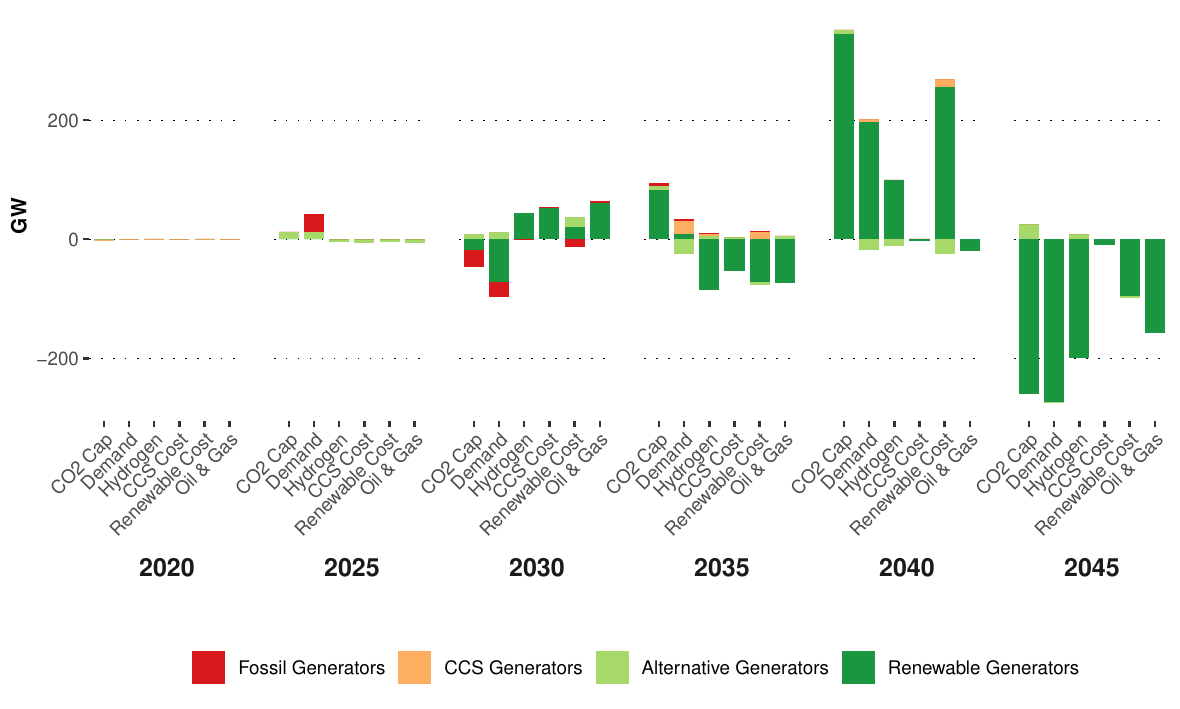}
    \caption{Structural difference in power generation investments between the deterministic solution and stochastic solution. (one long-term uncertainty)}
    \label{fig:One_dimension_diff}
\end{figure}

Moving on to the stochastic case, the expected investment decisions for each investment stage are compared with the deterministic solution in Figure \ref{fig:One_dimension_diff}. This illustrates the difference in the expected invested capacity for each stage. A positive value means larger expected investments for the MHSP than in the deterministic case, while a negative value means fewer expected investments. Note that the branching stages are in 2025, 2035 and 2045, meaning that the uncertainty is realised in 2030 and 2040. 

Little difference between the solutions is observed for the early stages before uncertainty is realised. For 2030, when the first scenario tree split is realised, the differences are more significant. A shift in investments between 2030 and 2035 is observed. For CO$_2$ cap and demand MHSPs, the shift results in more investments in 2035 than in 2030, while the opposite is the case for the other uncertain parameters. A similar effect can be seen in 2040 and 2045, where the differences are even more significant. The significant shift in renewables investments from 2045 to 2040 is consistent across six stochastic cases except for CCS Cost MHSP and Oil \& Gas MHSP. This shows how the structure of the scenario tree might affect the timing of the investments. 

Regarding the technology mix for the entire 30-year period, a general observation is that there is a small reduction in expected capacity investments in fossil generators and a slight increase in CCS generators when long-term uncertainty is modelled. 

\begin{figure}[!htb] 
    \centering
    \begin{subfigure}[t]{0.15\textwidth}
        \vspace{0.3cm} 
        \hspace{0.0cm} 
        \includegraphics[width=\textwidth]{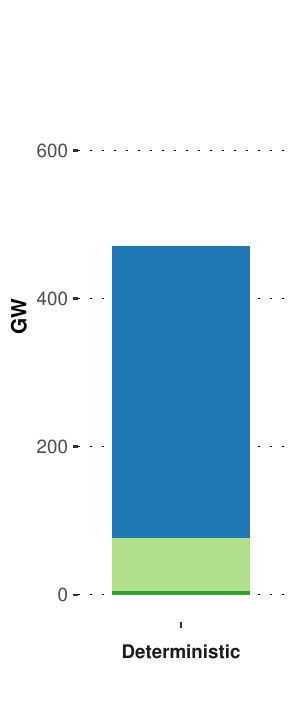}
        \caption{\scriptsize Deterministic}
        \label{fig:renewable_deterministic_2030}
    \end{subfigure}
    \hfill
    \begin{subfigure}[t]{0.80\textwidth}
        \vspace{0.2cm} 
        \includegraphics[width=\textwidth]{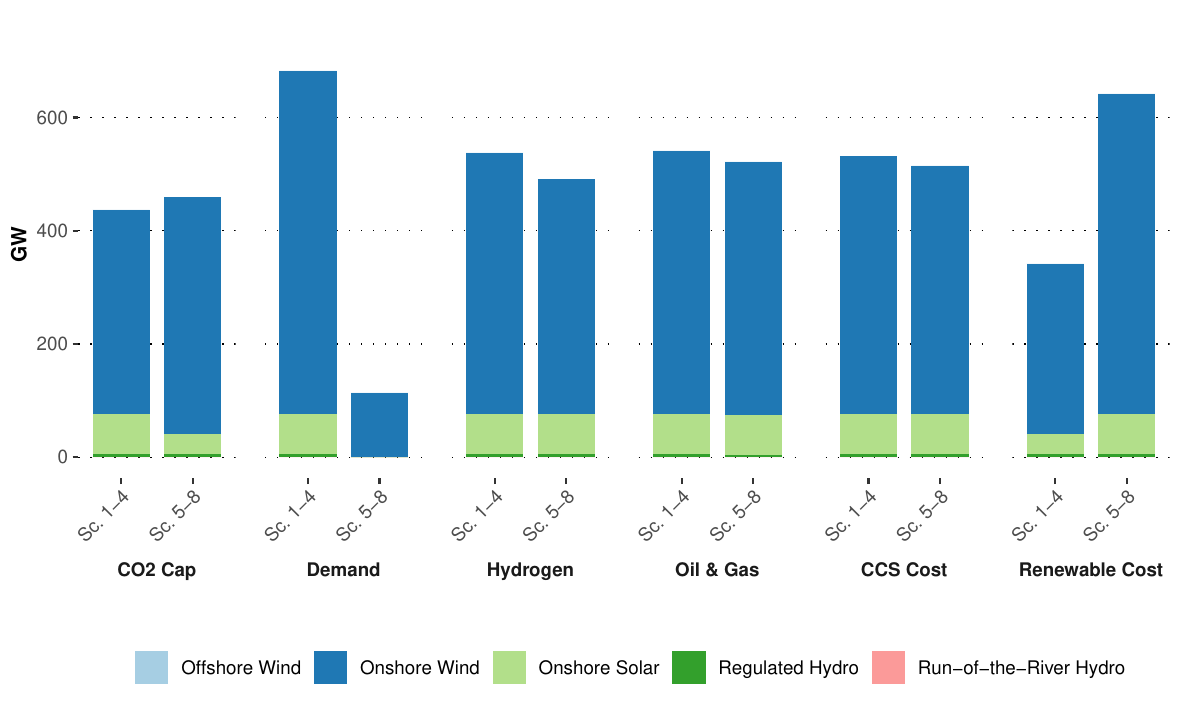}
        \caption{\scriptsize Stochastic}
        \label{fig:renewable_stochastic_2030}
    \end{subfigure}
    \caption{Investments in renewables in 2030 for the deterministic and different MHSP approaches (one long-term uncertainty).}
    \label{fig:investments_2030}
\end{figure}

The expected investments of the stochastic solutions offer insights into the impact of including uncertainty in the model. However, when real scenarios unfold, the solutions can diverge significantly from the expected outcomes. To gain a clearer insight into this, we analyse the scenarios at different branching stages. Figure \ref{fig:investments_2030} shows the structure of the investments in renewables in 2030. Every up-branch reflects a high realisation, while every down-branch reflects a low realisation of the uncertain parameter. If one considers the CO$_2$ Cap and Demand parameters, scenarios 1-4, the up-branch, reflects high allowed emission levels and high demand, while scenarios 5-8 represent low allowed emission levels and low demand. 



Figure \ref{fig:investments_2030} compares the investments with the deterministic solution. In 2030, onshore wind and some investments in onshore solar dominate the solution. For the models with uncertain CO$_2$ cap, hydrogen, oil \& gas and CCS cost, the structure of the solutions is nearly identical to the deterministic solutions, but with some more investments for both realisations of the uncertainty. However, the situation is different for the solutions with uncertain demand and renewable technology costs. As expected, the scenarios with high demand and the scenarios with lower renewable costs show increased investments.

\begin{figure}[!htb] 
    \centering
    \begin{subfigure}[t]{0.15\textwidth}
        \vspace{0.3cm} 
        \hspace{0.0cm} 
        \includegraphics[width=\textwidth]{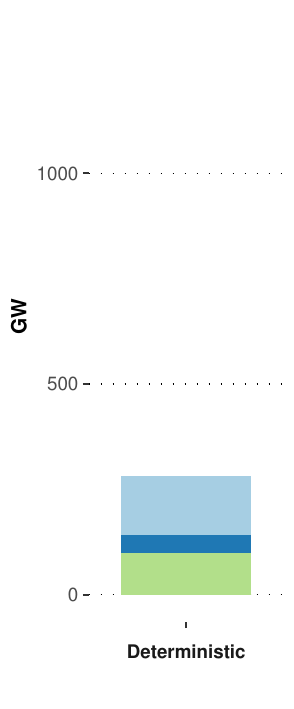}
        \caption{\scriptsize Deterministic}
        \label{fig:renewable_deterministic_2040}
    \end{subfigure}
    \hfill
    \begin{subfigure}[t]{0.80\textwidth}
        \vspace{0.2cm} 
        \includegraphics[width=\textwidth]{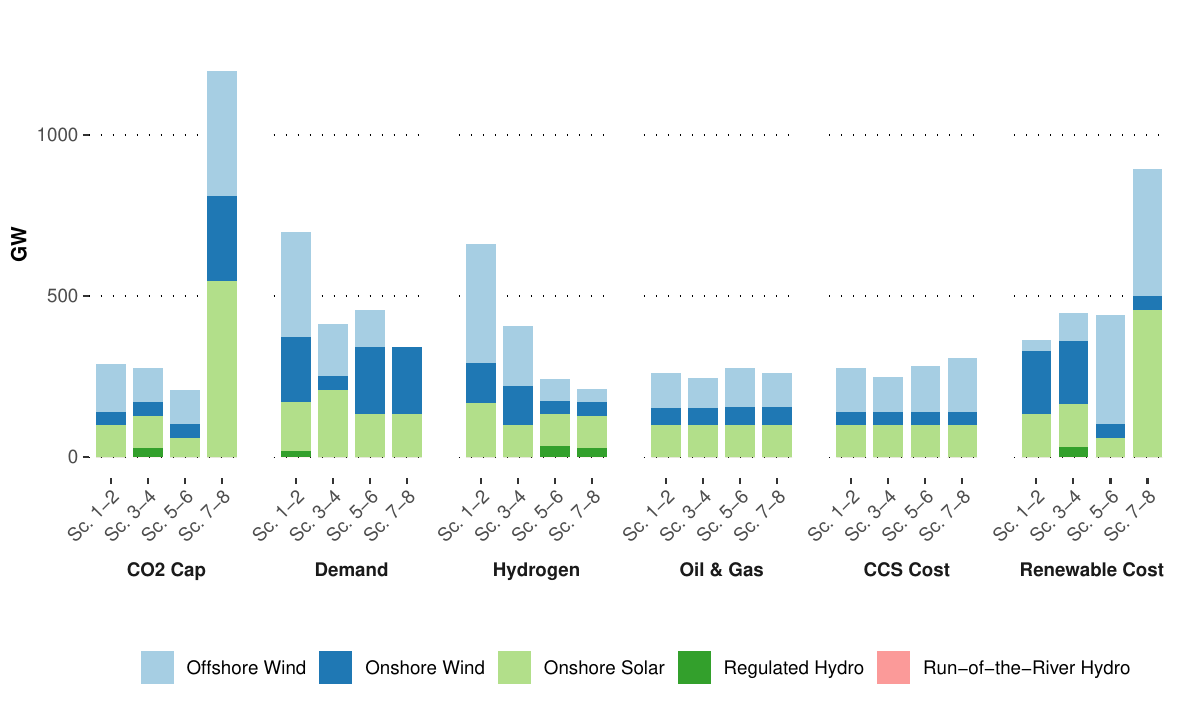}
        \caption{\scriptsize Stochastic}
        \label{fig:renewable_stochastic_2040}
    \end{subfigure}
    \caption{Investments in renewables in 2040 for the deterministic and different MHSP approaches.}
    \label{fig:investments_2040_renewables}
\end{figure}

For 2040, the differences between the branches for the individual SPs are more significant, as seen in Figure \ref{fig:investments_2040_renewables}. The differences between the individual scenarios are large for some of the uncertain parameters. For instance, the investments in renewables for the low-emission CO$_2$ scenarios are more than three times larger than for the other scenarios. This underlines the importance of assessing the scenarios individually, as these variations are not captured when looking at the expected solutions.

Another observation is that in the cases with uncertain demand and renewable costs, the investments are larger than that in the deterministic scenario for all realisations of uncertainty. In contrast, for scenarios with uncertain oil \& gas prices or CCS costs, the impact is less prominent, showing results very similar to the deterministic scenario across all scenario tree branches. This suggests that the connection between oil \& gas prices and CCS costs and renewable investments might not be as strong. Therefore, the influence of uncertainty in oil \& gas prices or CCS costs on investment decisions in renewables appears to be minimal compared to the effects of demand and renewable cost uncertainties. 

Lastly, some of the dynamics of the power system can be interpreted from the hydrogen MHSP. It is clear that in scenarios with high hydrogen demand, there are also increased investments in renewable energy, which might be connected with an increased need for electrolysis, which has good synergies with intermittent renewable energy. 

\begin{figure}[!htb]
    \centering
    \includegraphics[height=6cm]{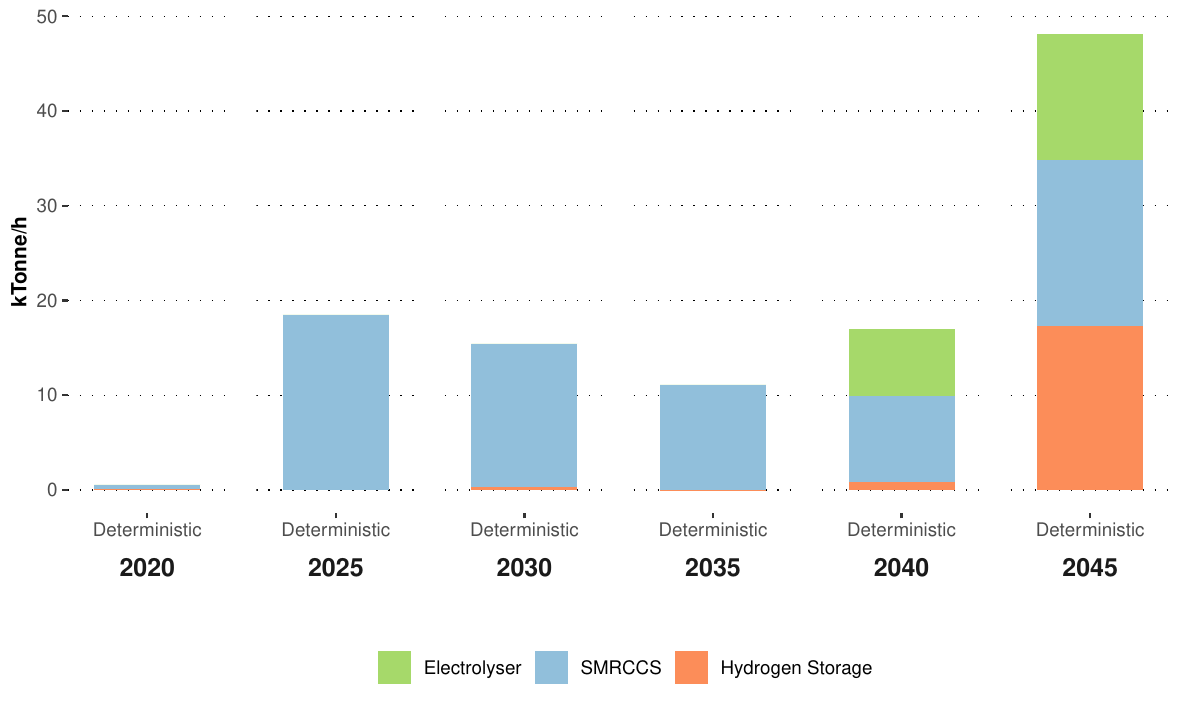}
    \caption{Investments in hydrogen generation technologies and storage.}
    \label{fig:Deterministic_hydrogen_generation}
\end{figure}

Figure \ref{fig:Deterministic_hydrogen_generation} shows the evolution of investments in hydrogen generation and storage technologies in the deterministic case. Note here that the units are different for the two, with hydrogen storage being measured in total storage capacity, while production technologies are measured in hourly production capacity. Nevertheless, it is interesting to compare the two as it makes it possible to compare the total storage capacity with the hourly production capacity. In the deterministic model, SMRCCS dominates the solution, especially in the early stages before investments in electrolysis grow towards the end of the horizon. Hydrogen storage is nearly negligible until 2045 when some investments are undertaken. Compared with the total production capacity, the storage capacity is relatively low.

\begin{figure}[!htb]
    \centering
    \includegraphics[width=\textwidth]{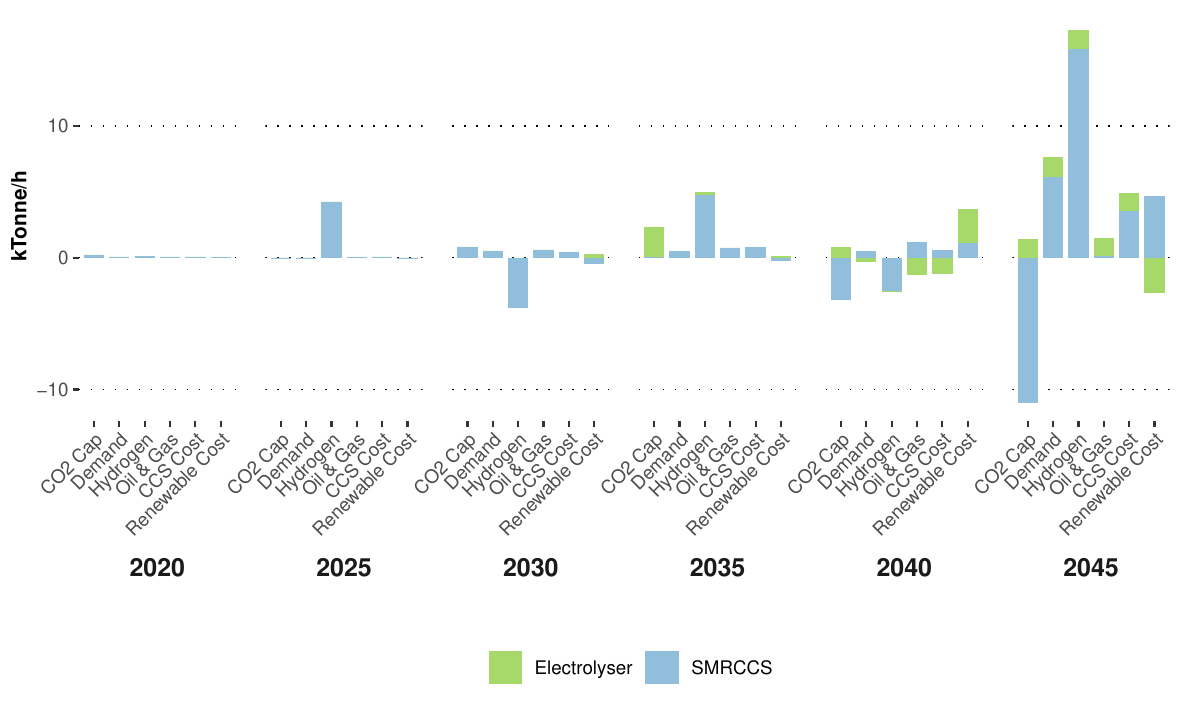}
    \caption{Structural difference in hydrogen generation investments between the deterministic solution and one-dimensional stochastic solution.}
    \label{fig:One_dimension_hydrogen_generation}
\end{figure}

In Figure \ref{fig:One_dimension_hydrogen_generation}, the structural differences between the deterministic and expected MHSP solution for hydrogen production are presented. The scenario solutions generally align closely with the deterministic solution during the first four periods. However, there is an exception for the scenario involving stochastic hydrogen demand, which shows a shift in investments in SMRCSS from 2030 to 2025 and from 2040 to 2035. In 2045, the hydrogen related to a scenario tree expects significantly more investments than the deterministic instance. The opposite is the case for the CO$_2$ emission constraint scenario, which results in less investment in SMRCCS while slightly increasing investments in electrolysis, without making up for the reduced SMRCCS investments.

\begin{figure}[!htb]
    \centering
    \includegraphics[width=\textwidth]{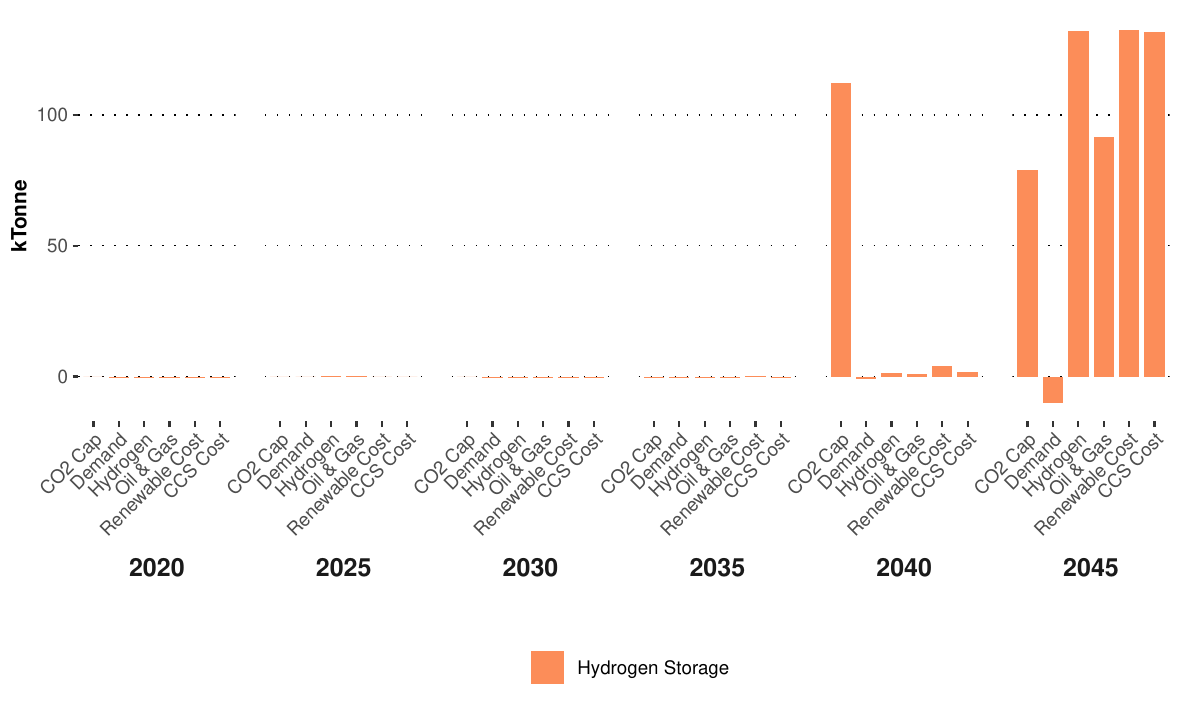}
    \caption{Structural difference in hydrogen storage investments between the deterministic solution and one-dimensional stochastic solution.}
    \label{fig:One_dimension_hydrogen_storage}
\end{figure}

Analysing the disparity between the deterministic and stochastic solutions for hydrogen storage, as depicted in Figure \ref{fig:One_dimension_hydrogen_storage}, we observe practically no change in the first four investment stages. However, a significant shift occurs towards the end of the horizon, with substantial investments in storage across most MHSP solutions, with the exception of demand uncertainty. A potential explanation for this pattern is a growing need for enhanced flexibility in response to increased uncertainty and more extreme outcomes. Hydrogen is a promising solution in this context, given its versatility as an energy carrier, and this can be considered a hedging strategy for the system. Especially interesting might be the CO$_2$ cap scenario, which reduced overall hydrogen production capacity. This is in high contrast to hydrogen storage investments that increase significantly. 

\begin{figure}[!htb] 
    \centering
    \begin{subfigure}[t]{0.15\textwidth}
        \vspace{0.3cm} 
        \hspace{0.0cm} 
        \includegraphics[width=\textwidth]{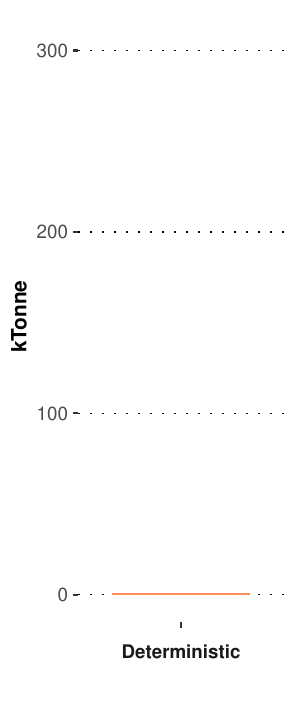}
        \caption{\scriptsize Deterministic}
    \end{subfigure}
    \hfill
    \begin{subfigure}[t]{0.80\textwidth}
        \vspace{0.2cm} 
        \includegraphics[width=\textwidth]{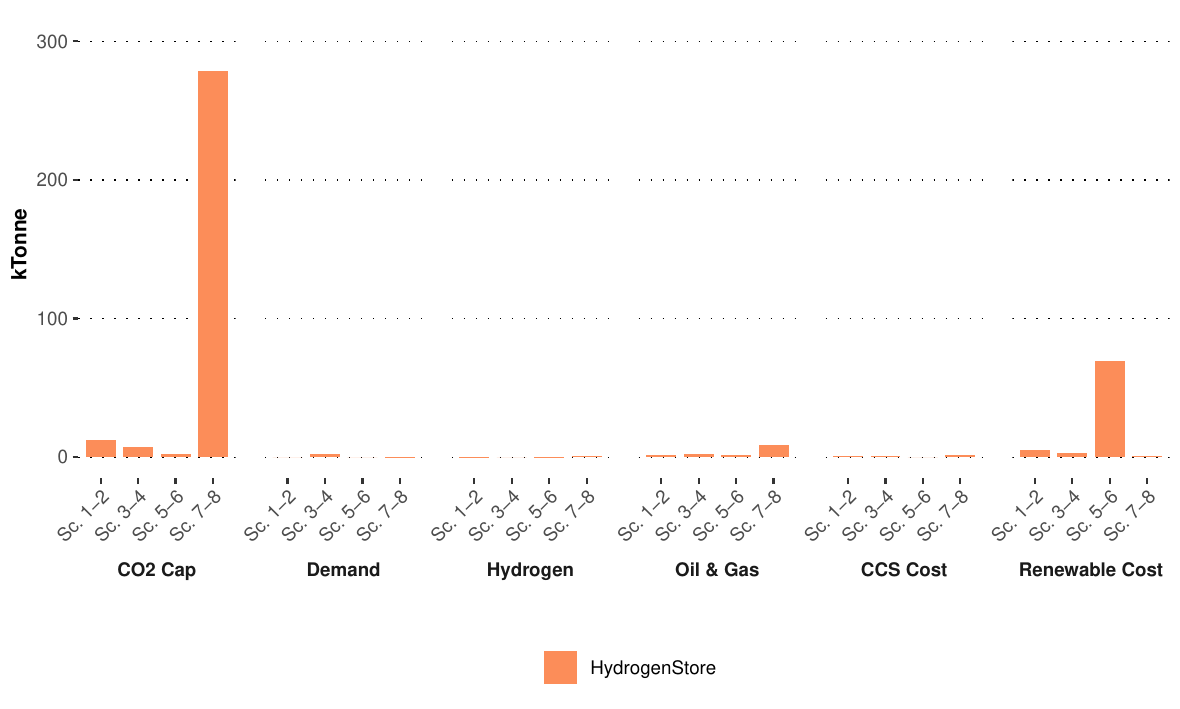}
        \caption{\scriptsize Stochastic}
    \end{subfigure}
    \caption{Investments in hydrogen storage in 2040 for the deterministic and different MHSP approaches.}
    \label{fig:investments_2040_hydrogen_storage}
\end{figure}

Diving deeper, we can see that for CO$_2$ budget uncertainty, the large investment in hydrogen storage in 2040 fully coincides with the instance with drastically higher investments in renewable energy sources from Figure \ref{fig:investments_2040_hydrogen_storage}. Further examination reveals the same trend for 2045, where the CO$_2$ cap instance with high investments in renewable energy also had high investments in energy storage. An interpretation of these factors, well backed up by real-world dynamics, is that high renewable energy penetration causes increased utilisation of electrolysers, which again causes a need for balancing hydrogen load through storage investments. This again showcases the value of stochastic planning, as tail scenarios and the consequences of these can be captured by the planning model. Still, similar effects are not as clear for other uncertain parameters, such as renewable costs, which do not see the same clear connection between high renewable penetration and high hydrogen storage investments.  

\subsection{Combined long-term uncertainty}\label{two_dimensional_uncertainty}
This section expands the analysis of combined uncertainty. We consider the CO$_2$ cap and demand uncertainty simultaneously. The scenario trees used in this assessment assume no correlation between the parameters, consequently increasing the number of scenarios from 8 to 64. The parameter values for the scenarios are held at the same level as for the one-dimensional uncertainty.

\begin{figure}[!htb]
    \centering
    \includegraphics[width=\textwidth]{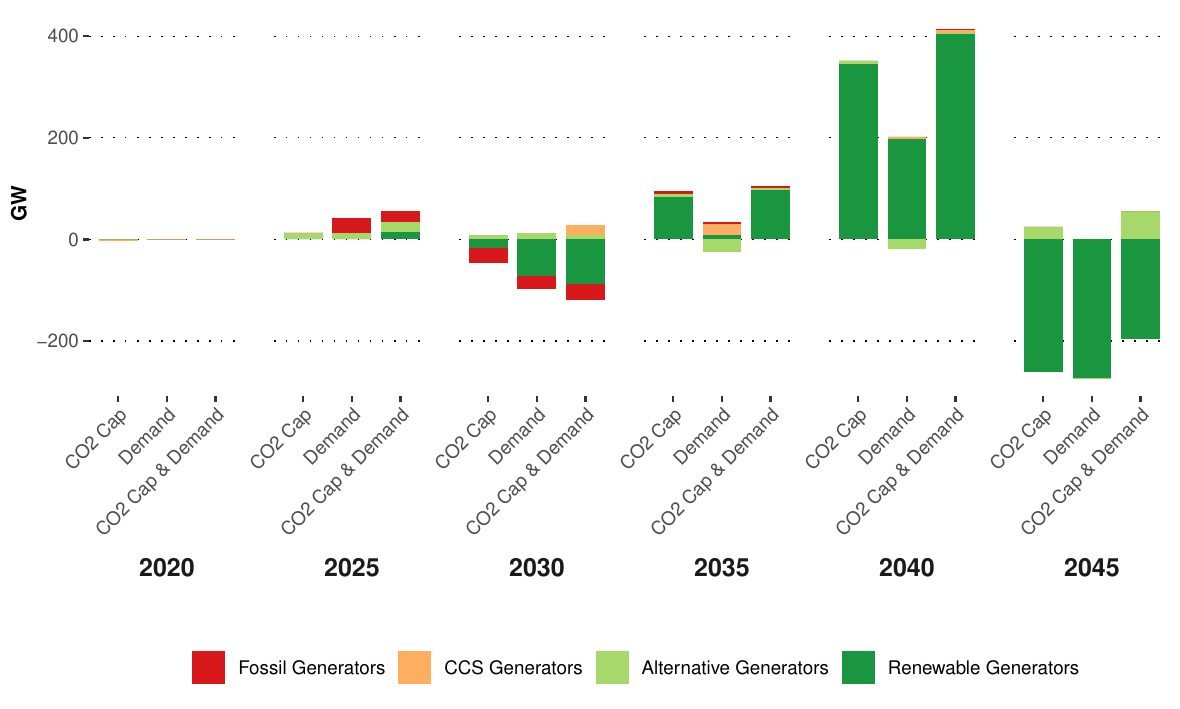}
    \caption{Investments in power generation for combined and individual uncertainty of the CO$_2$ cap and demand scenarios.}
    \label{fig:2d_co2_demand}
\end{figure}

Figure \ref{fig:2d_co2_demand} illustrates the difference in power generation investments between the deterministic scenario and scenarios with uncertain parameters. It compares the results with the combined uncertainty of CO$_2$ cap and demand with the individual uncertainties of each parameter.

The above figures and analysis indicate that the combined uncertainties do not give very different results than individual assessments. Few significant differences are observed when looking at the quantities and types of investments. This trend is observed across all combined uncertainty settings. 



\subsection{Assessment of long-uncertainty using VSS and RHVSS}\label{RH_section}
In this section, we show the value of long-term uncertainty modelling by assessing VSS and RHVSS.  The calculation is made with uncertainty on the CO$_2$ cap.

When VSS is used, the objective of the expected value solution is 448.85 €bn compared with 346.16 €bn from the MHSP solution, which yields 102.69 €bn VSS, around 29.7\%.

For RHVSS, Table \ref{Objective_rh_co2} shows the objective for each of the eight scenarios with CO$_2$ uncertainty, solved using the RH approach. Scenario 1 is the scenario for high-emission, whereas scenario 8 is the low-emission scenario. The scenarios are weighted differently. The RHVSS is around  1.7\%. This shows that the value of the MHSP solution compared to the RH approach is limited despite giving a positive value. It is worth mentioning here that the RH approach, despite its simplicity, is much less naive than solving the deterministic EEV problem. 

\begin{table}[!htb]
\centering
\caption{Objective value for the different RH solutions, with CO$_2$ cap uncertainty.}
\footnotesize
\begin{tabular}{c|S[table-format=3.2]|S[table-format=3.2]|S[table-format=3.2]|S[table-format=1.2]}
\toprule
{Scenario} & {Objective (€bn)} & {Investment cost (€bn)} & {Operational cost (€bn)} & {Weight} \\ \hline
1 & 320.82 & 157.26 & 163.56 & 0.05 \\
2 & 320.76 & 157.30 & 163.46 & 0.05 \\
3 & 322.27 & 159.76 & 162.51 & 0.10 \\
4 & 326.51 & 156.85 & 169.66 & 0.10 \\
5 & 328.44 & 182.21 & 146.23 & 0.10 \\
6 & 341.73 & 182.22 & 159.51 & 0.20 \\
7 & 381.14 & 200.46 & 180.68 & 0.20 \\
8 & 389.32 & 200.70 & 188.62 & 0.20 \\ \hline
{ERHEV} & 352.24 & 182.29 & 169.95 & {-} \\ \hline
{Stochastic} & 346.16 & 179.13 & 167.03 & {-} \\ \hline
{RHVSS} & 6.08 & {-} & {-} & {-}\\
\bottomrule
\end{tabular}
\label{Objective_rh_co2}
\end{table}

\begin{figure}[!htb]
    \centering
    \includegraphics[width=\textwidth]{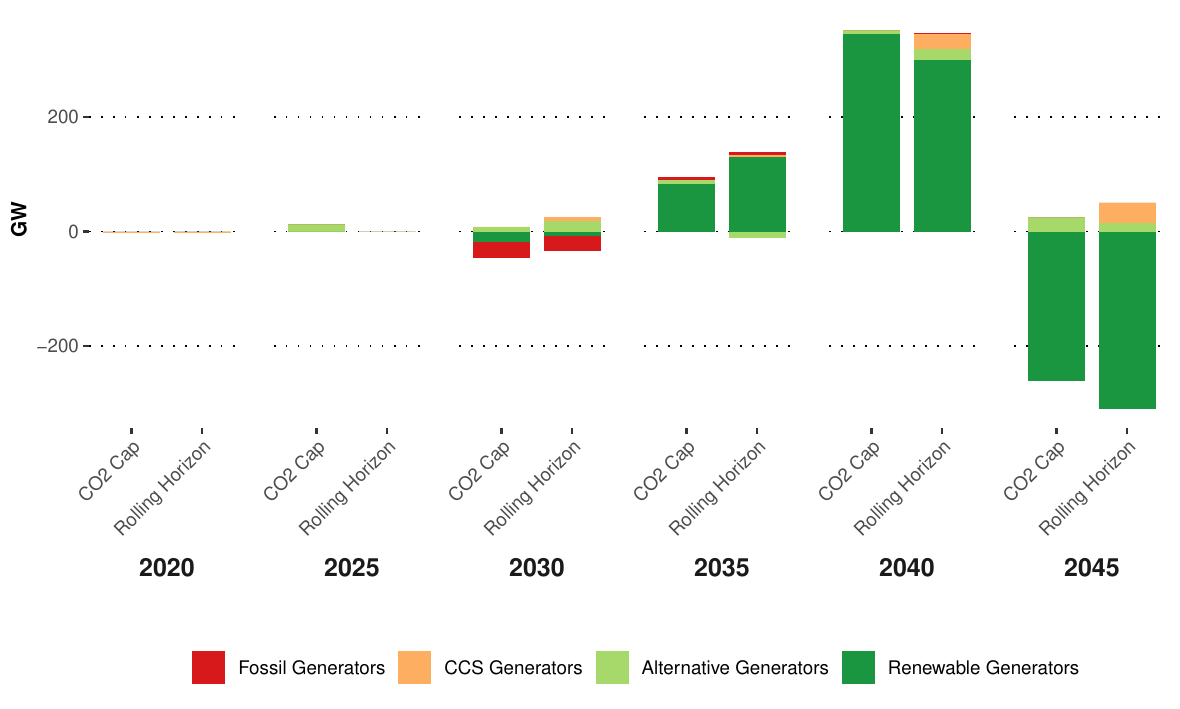}
    \caption{Investments in power generation for the MHSP and RH approach, compared with the deterministic baseline.}
    \label{fig:diff_rolling_horizon}
\end{figure}

To get a better insight into the differences between the MHSP and RH approaches, we compare the investment decisions. Figure \ref{fig:diff_rolling_horizon} shows a plot of the investments in power generation compared with the deterministic baseline. As can be seen, the investments in renewables are relatively similar for both approaches. Another observation is that the RH approach invests more in CCS technologies and alternative generators than the MHSP. 

Both methods show similar behaviour regarding when decisions are made, especially for the last two periods when investments shift from 2045 to 2040. This is a natural consequence of when the information in the scenario tree is revealed. Both methods receive information about which scenario they belong to and update the future expected values simultaneously.    

\begin{figure}[!htb] 
    \centering
    \begin{subfigure}[t]{0.223\textwidth}
        \vspace{0.2cm} 
        \hspace{0.0cm} 
        \includegraphics[width=\textwidth]{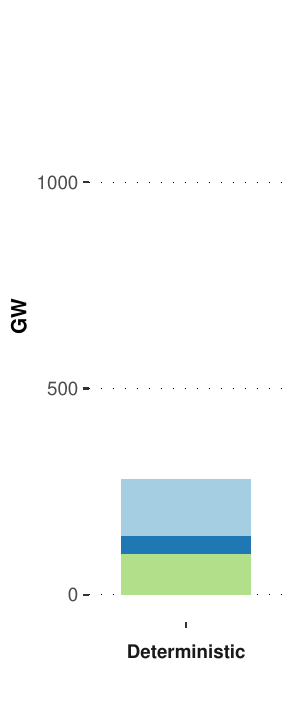}
        \caption{\scriptsize Deterministic}
        \label{fig:renewable_deterministic_RH_2040}
    \end{subfigure}
    \hfill
    \begin{subfigure}[t]{0.75\textwidth}
        \vspace{0.2cm} 
        \includegraphics[width=\textwidth]{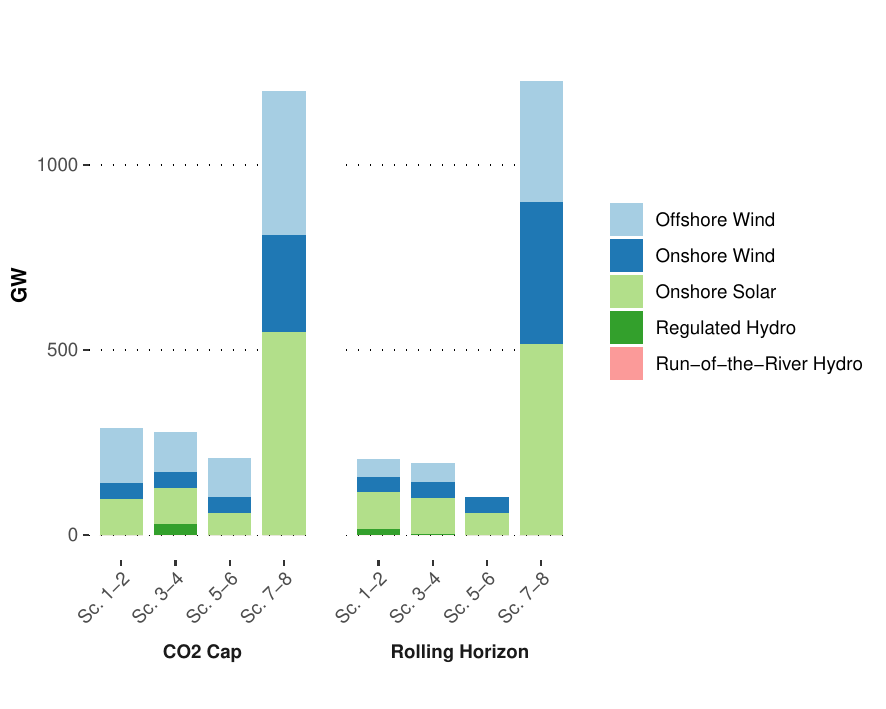}
        \caption{\scriptsize Stochastic}
        \label{fig:renewable_RH_2040}
    \end{subfigure}
    \caption{Comparison of renewable investments in 2040 for the RH and MHSP approach with uncertain CO$_2$ constraint.}
    \label{fig:investments_RH_2040}
\end{figure}

A comparison of the scenario-specific investments in renewables between the RH and MHSP approach with uncertain CO$_2$ cap constraint is made in Figure \ref{fig:investments_RH_2040}. The structures of the two solutions are similar, both when comparing the size and type of investments. The low-emission scenario requires significantly larger investments in both cases. Onshore and offshore wind, and onshore solar dominate the investments in all MHSP and RH solutions.  

\section{Conclusions and future work}
\label{sec:conclusions}
This paper comprehensively modelled and analysed the impact of multi-timescale uncertainty in energy system planning problems. We extended the REORIENT model for the analysis by introducing random sampling and moment matching SGRs for the short-term scenarios and a long-term SGR. The REORIENT model is a multi-horizon stochastic MILP for integrated investment, retrofit and abandonment energy system planning. To address the computational difficulty, we proposed parallel stabilised Benders decomposition. We applied the proposed model and algorithm to the European energy system planning problem. We considered uncertainty in short-term time series parameters, and long-term uncertainty parameters, including CO$_2$ cap, power demand, hydrogen demand, oil and gas prices, CCS technology cost, and renewable technology cost. We evaluated the stability of the stochastic solutions, VSS and RHVSS, objectives and system investment decisions. The results showed that: (1)including multi-timescale uncertainty yields a significant value of the stochastic solutions, (2) long-term uncertainty in the right-hand side parameters affects the solution structure than cost coefficient uncertainty, (3) parallel stabilised Benders decomposition is up to 7.5 times faster than the serial version.

Generating scenarios that accurately represent long-term uncertainty is challenging. Therefore, a possible future work would be to develop and test different scenario generation methods for long-term uncertainty. 

\section*{CRediT author statement}
\textbf{Hongyu Zhang:} Conceptualisation, Methodology, Software, Validation, Formal analysis, Investigation, Visualisation, Data curation, Writing - original draft, Writing - review, \& editing, Supervision. \textbf{Asbjørn Nisi:} Conceptualisation, Methodology, Software, Validation, Formal analysis, Investigation, Visualisation, Data curation. \textbf{Erlend Heir:} Conceptualisation, Methodology, Software, Validation, Formal analysis, Investigation, Visualisation, Data curation. \textbf{Asgeir Tomasgard:} Conceptualisation, Methodology, Supervision, Funding acquisition.

\section*{Declaration of competing interest}
The authors declare that they have no known competing financial interests or personal relationships that could have appeared to influence the work reported in this paper.






\clearpage

\setlength{\bibsep}{0pt plus 0.3ex}
\footnotesize{
\bibliographystyle{model5-names}
\bibliography{multi-timescale_uncertainty_energy_system_planning_arxiv}}

\end{document}